\documentclass[letterpaper,american,english]{amsart}
\usepackage[T1]{fontenc}
\usepackage[latin1]{inputenc}
\usepackage{fancyhdr}
\usepackage{amsthm}
\usepackage{amstext}
\usepackage{amssymb}
\usepackage[all]{xy}

\makeatletter
\numberwithin{equation}{section}
\numberwithin{figure}{section}
\theoremstyle{plain}
\newtheorem{thm}{Theorem}[section]
  \theoremstyle{definition}
  \newtheorem{defn}[thm]{Definition}
  \theoremstyle{remark}
  \newtheorem*{rem*}{Remark}
  \theoremstyle{remark}
  \newtheorem{rem}[thm]{Remark}
 \theoremstyle{definition}
  \newtheorem{example}[thm]{Example}
  \theoremstyle{plain}
  \newtheorem{prop}[thm]{Proposition}
  \theoremstyle{plain}
  \newtheorem{lem}[thm]{Lemma}
  \theoremstyle{definition}
  \newtheorem{condition}[thm]{Condition}
  \theoremstyle{plain}
  \newtheorem{cor}[thm]{Corollary}

\usepackage{verbatim} 
\usepackage{eucal}
\usepackage{amsthm}
\usepackage[all]{xy}
\usepackage{manfnt}
\usepackage{float}
\usepackage{xr}
\SelectTips{cm}{}
\xyoption{line}
\usepackage{lmodern}
\usepackage{hyperref}
\makeatletter \newcommand{\xyR}[1]{%
\makeatletter \xydef@\xymatrixrowsep@{#1} \makeatother }
\makeatletter \newcommand{\xyC}[1]{%
\makeatletter \xydef@\xymatrixcolsep@{#1} \makeatother }
\mathcode`\:="603A 
\DeclareSymbolFont{rsfs}{U}{rsfs}{m}{n}
\DeclareSymbolFontAlphabet{\mathrf}{rsfs}

\makeatother

\begin{document}

\title{Model-categories of coalgebras over operads}

\author{Justin R. Smith}



\curraddr{Department of Mathematics\\
Drexel University\\
Philadelphia,~PA 19104}

\email{jsmith@drexel.edu}

\urladdr{http://vorpal.math.drexel.edu}

\date{\today}
\begin{abstract}
This paper constructs model structures on the categories of coalgebras
and pointed irreducible coalgebras over an operad whose components
are projective, finitely generated in each dimension, and satisfy
a condition that allows one to take tensor products with a unit interval.
The underlying chain-complex is assumed to be unbounded and the results
for bounded coalgebras over an operad are derived from the unbounded
case.
\end{abstract}
\maketitle
\tableofcontents{}

\global\long\def\ring{R}
\global\long\def\integers{\mathbb{Z}}
\global\long\def\betabar{\bar{\beta}}
 \global\long\def\desusp{\downarrow}
\global\long\def\susp{\uparrow}
\global\long\def\cobar{\mathcal{F}}
\global\long\def\coend{\mathrm{CoEnd}}
\global\long\def\ainfty{A_{\infty}}
\global\long\def\coassoc{\mathrm{Coassoc}}
\global\long\def\trm{\mathrm{T}}
\global\long\def\tfr{\mathfrak{T}}
\global\long\def\tabbr{\hat{\trm}}
\global\long\def\Tabbr{\hat{\tfr}}
\global\long\def\afr{\mathfrak{A}}
\global\long\def\homz{\mathrm{Hom}_{\ring}}
\global\long\def\zend{\mathrm{End}}
\global\long\def\rs#1{\mathrm{R}S_{#1 }}
\global\long\def\forgetful#1{\lceil#1\rceil}
\global\long\def\highprod#1{\bar{\mu}_{#1 }}
\global\long\def\slength#1{|#1 |}
\global\long\def\barcs{\bar{\mathcal{B}}}
\global\long\def\ubarcs{\mathcal{B}}
\global\long\def\zs#1{\ring S_{#1 }}
\global\long\def\homzs#1{\mathrm{Hom}_{\ring S_{#1 }}}
\global\long\def\zpi{\mathbb{Z}\pi}
\global\long\def\D{\mathfrak{D}}
\global\long\def\ahat{\hat{\mathfrak{A}}}
\global\long\def\cbar{{\bar{C}}}
\global\long\def\cf#1{\mathcal{C}(#1 )}
\global\long\def\ddelta{\dot{\Delta}}
\global\long\def\dimlimiter{\triangleright}
\global\long\def\coalgcat{\mathrf S_{0}}
\global\long\def\hcoalgcat{\mathrf{S}}
\global\long\def\ircoalgcat{\mathrf I_{0}}
\global\long\def\bircoalgcat{\mathrf{I}_{0}^{+}}
\global\long\def\hircoalgcat{\mathrf I}
\global\long\def\dcoalgcat{\mathrm{ind}-\coalgcat}
\global\long\def\chaincat{\mathbf{Ch}}
\global\long\def\coll{\mathrm{Coll}}
\global\long\def\bchaincat{\mathbf{Ch}_{0}}
\global\long\def\ilimit{\varprojlim\,}
\global\long\def\bigboxtimes{\mathop{\boxtimes}}
\global\long\def\dlimit{\varinjlim\,}
\global\long\def\coker{\mathrm{{coker}}}
\global\long\def\icoalgcat{\mathrm{pro}-\mathrf{S}_{0}}
\global\long\def\iircoalgcat{\mathrm{pro-}\ircoalgcat}
\global\long\def\dircoalgcat{\mathrm{ind-}\ircoalgcat}
\global\long\def\core#1{\left\langle #1\right\rangle }
\global\long\def\ilimitder{\varprojlim^{1}\,}
\global\long\def\pcoalg#1#2{P_{\mathcal{#1}}(#2) }
\global\long\def\pcoalgf#1#2{P_{\mathcal{#1}}(\forgetful{#2}) }
\global\long\def\coequalizer{\mathop{\mathrm{coequalizer}}}
\global\long\def\mainoperad{\mathcal{H}}
\global\long\def\cone#1{\mathrm{Cone}(#1)}
\global\long\def\im{\operatorname{im}}

\section{Introduction}

\newdir{ >}{{}*!/-5pt/@{>}}

\lhead{HOMOTOPY THEORY OF COALGEBRAS}

\rhead{JUSTIN R. SMITH}

Although the literature contains several papers on homotopy theories
for \emph{algebras} over operads --- see~\cite{Kriz-May}, \cite{Mandell-einfty},
and \cite{operad-book} --- it is more sparse when one pursues similar
results for coalgebras. In \cite{Quillen:1969}, Quillen developed
a model structure on the category of 2-connected cocommutative coalgebras
over the rational numbers. V.~Hinich extended this in \cite{Hinich-DG-coalgebras}
to coalgebras whose underlying chain-complexes were unbounded (i.e.,
extended into negative dimensions). Expanding on Hinich's methods,
K.~Lefèvre derived a model structure on the category of coassociative
coalgebras --- see \cite{Coalgebra-thesis}. In general, these authors
use indirect methods, relating coalgebra categories to other categories
with known model structures.

Our paper finds model structures for coalgebras over any operad fulfilling
a basic requirement (condition~\ref{cond:mainoperadassumption}).
Since operads uniformly encode many diverse coalgebra structures (coassociative-,
Lie-, Gerstenhaber-coalgebras, etc.), our results have wide applicability. 

The author's intended application involves investigating the extent
to which Quillen's results in rational homotopy theory (\cite{Quillen:1969})
can be generalized to \emph{integral} homotopy theory.

Several unique problems arise that require special techniques. For
instance, constructing injective resolutions of coalgebras naturally
leads into infinitely many negative dimensions. The resulting model
structure --- and even that on the underlying chain-complexes ---
fails to be cofibrantly generated (see \cite{Quillen-model-rel-homolog-alg}).
Consequently, we cannot easily use it to induce a model structure
on the category of coalgebras.

We develop the general theory for unbounded coalgebras, and derive
the bounded results by applying a truncation functor.

In \S~\ref{sec:Notation-and-conventions}, we define operads and
coalgebras over operads. We also give a basic condition (see \ref{cond:mainoperadassumption})
on the operad under consideration that we assume to hold throughout
the paper. Cofibrant operads always satisfy this condition and every
operad is weakly equivalent to one that satisfies this condition.

In \S~\ref{sec:Model-categories}, we briefly recall the notion
of model structure on a category and give an example of a model structure
on the category of unbounded chain-complexes. 

In $\S$~\ref{sec:The-general-case}, we define a model structure
on categories of coalgebras over operads. When the operad is projective
and finitely-generated in all dimensions, we verify that nearly free
coalgebras satisfy Quillen's axioms of a model structure (see \cite{Quillen:1967}
or \cite{Goerss-Jardine}).

Section~\ref{sub:Description-of-the} describes our model-structure
--- classes of cofibrations, fibrations and weak equivalences. Section~\ref{sub:Proof-of-CM1}
proves the first few axioms of a model-structure (CM~1 through CM~3,
in Quillen's notation). Section~\ref{sub:Proof-of-CM5} proves axiom
CM~5, and section~\ref{sub:Proof-of-CM4} proves CM~4.

A key step involves proving the existence of cofibrant and fibrant
replacements for objects. In our model structure, all coalgebras are
cofibrant (solving this half of the problem) and the hard part of
is to find \emph{fibrant} replacements. 

We develop resolutions of coalgebras by cofree coalgebras --- our
so-called rug-resolutions --- that solves the problem: see lemma~\ref{lem:cofreeresolution}
and corollary~\ref{cor:factortrivialcofibrationfibration}. This
construction naturally leads into infinitely many negative dimensions
and was the motivation for assuming underlying chain-complexes are
unbounded. 

All coalgebras are cofibrant and \emph{fibrant} coalgebras are characterized
as retracts of canonical resolutions called rug-resolutions (see corollary~\ref{cor:fibrationcharacterization}
and corollary~\ref{cor:fibrationcharacterization}) --- an analogue
to total spaces of Postnikov towers.

In the cocommutative case over the rational numbers, the model structure
that we get is \emph{not} equivalent to that of Hinich in \cite{Hinich-DG-coalgebras}.
He gives an example (9.1.2) of a coalgebra that is acyclic but not
contractible. In our theory it \emph{would} be contractible, since
it is over the rational numbers and bounded.

In \S~\ref{sec:The-bounded-case}, we discuss the (minor) changes
to the methods in \S~\ref{sec:The-general-case} to handle coalgebras
that are bounded from below. This involves replacing the cofree coalgebras
by their truncated versions.

In \S~\ref{sec:Examples}, we consider two examples over the rational
numbers. In the rational, 2-connected, cocommutative, coassociative
case, we recover the model structure Quillen defined in \cite{Quillen:1969}
--- see example~\ref{exa:quillen}.

In appendix~\ref{sec:nearlyfree}, we study \emph{nearly free} $\integers$-modules.
These are modules whose countable submodules are all $\integers$-free.
They take the place of free modules in our work, since the cofree
coalgebra on a free modules is not free (but is nearly free).

In appendix~\ref{sec:Basic-constructions}, we develop essential
category-theoretic constructions, including equalizers (\S~\ref{sub:Core-of-a}),
products and fibered products (\S~\ref{sub:Categorical-products}),
and colimits and limits (\S~\ref{sub:Limits-and-colimits}). The
construction of limits in \S~\ref{sub:Limits-and-colimits} was
this project's most challenging aspect and consumed the bulk of the
time spent on it. This section's key results are corollary~\ref{cor:intersectionalreadycoalg},
which allows computation of inverse limits of coalgebras and theorem~\ref{thm:aequivforgetfuldirectsummand},
which shows that these inverse limits share a basic property with
those of chain-complexes.

I am indebted to Professor Bernard Keller for several useful discussions.

\section{Notation and conventions\label{sec:Notation-and-conventions}}

\label{sec:generalconstruction}

Throughout this paper, $\ring$ will denote a field or $\integers$. 
\begin{defn}
\label{def:nearlyfree}An $\ring$-module $M$ will be called \emph{nearly
free} if every countable submodule is $\ring$-free.\end{defn}
\begin{rem*}
This condition is automatically satisfied unless $\ring=\integers$.

Clearly, any $\integers$-free module is also nearly free. The Baer-Specker
group, $\integers^{\aleph_{0}}$, is a well-known example of a nearly
free $\integers$-module that is \emph{not} free --- see \cite{Fuchs:1970},
\cite{Baer:1937}, and \cite{Baer-Specker-nonfree}. Compare this
with the notion of $\aleph_{1}$-\emph{free} \emph{groups} --- see
\cite{Blass-Gobel:1996}. 

By abuse of notation, we will often call chain-complexes nearly free
if their underlying modules are (ignoring grading).

Nearly free $\integers$-modules enjoy useful properties that free
modules do \emph{not}. For instance, in many interesting cases, the
cofree coalgebra of a nearly free chain-complex is nearly free.\end{rem*}
\begin{defn}
\label{def:chaincat}We will denote the closed symmetric monoidal
category of unbounded, nearly free $\ring$-chain-complexes with $\ring$-tensor
products by $\chaincat$. We will denote the category of $\ring$-\emph{free}
chain chain-complexes that are\emph{ bounded from below} in dimension
$0$ by $\bchaincat$.
\end{defn}
The chain-complexes of $\chaincat$ are allowed to extend into arbitrarily
many negative dimensions and have underlying graded $\ring$-modules
that are
\begin{itemize}
\item arbitrary if $\ring$ is a field (but they will be free)
\item \emph{nearly free,} in the sense of definition~\ref{def:nearlyfree},
if $\ring=\integers$.
\end{itemize}
We make extensive use of the Koszul Convention (see~\cite{Gugenheim:1960})
regarding signs in homological calculations:
\begin{defn}
\label{def:koszul} If $f:C_{1}\to D_{1}$, $g:C_{2}\to D_{2}$ are
maps, and $a\otimes b\in C_{1}\otimes C_{2}$ (where $a$ is a homogeneous
element), then $(f\otimes g)(a\otimes b)$ is defined to be $(-1)^{\deg(g)\cdot\deg(a)}f(a)\otimes g(b)$. \end{defn}
\begin{rem}
If $f_{i}$, $g_{i}$ are maps, it isn't hard to verify that the Koszul
convention implies that $(f_{1}\otimes g_{1})\circ(f_{2}\otimes g_{2})=(-1)^{\deg(f_{2})\cdot\deg(g_{1})}(f_{1}\circ f_{2}\otimes g_{1}\circ g_{2})$.\end{rem}
\begin{defn}
\label{def:unitinterval}The symbol $I$ will denote the unit interval,
a chain-complex given by\begin{eqnarray*}
I_{0} & = & \ring\cdot p_{0}\oplus\ring\cdot p_{1}\\
I_{1} & = & \ring\cdot q\\
I_{k} & = & 0\,\text{if }k\ne0,1\\
\partial q & = & p_{1}-p_{0}\end{eqnarray*}
Given $A\in\chaincat$, we can define\[
A\otimes I\]
and \[
\cone A=A\otimes I/A\otimes p_{1}\]

\end{defn}
The set of morphisms of chain-complexes is itself a chain complex:
\begin{defn}
\label{def:homcomplex}Given chain-complexes $A,B\in\chaincat$ define\[
\homz(A,B)\]
to be the chain-complex of graded $\ring$-morphisms where the degree
of an element $x\in\homz(A,B)$ is its degree as a map and with differential\[
\partial f=f\circ\partial_{A}-(-1)^{\deg f}\partial_{B}\circ f\]
As a $\ring$-module $\homz(A,B)_{k}=\prod_{j}\homz(A_{j},B_{j+k})$.\end{defn}
\begin{rem*}
Given $A,B\in\chaincat^{S_{n}}$, we can define $\homzs n(A,B)$ in
a corresponding way.\end{rem*}
\begin{defn}
\label{def:operad}If $G$ is a discrete group, let $\bchaincat^{G}$
denote the category of chain-complexes equipped with a right $G$-action.
This is again a closed symmetric monoidal category and the forgetful
functor $\bchaincat^{G}\to\bchaincat$ has a left adjoint, $(-)[G]$.
This applies to the symmetric groups, $S_{n}$, where we regard $S_{1}$
and $S_{0}$ as the trivial group. The \emph{category of collections}
is defined to be the product\[
\mathrm{Coll}(\bchaincat)=\prod_{n\ge0}\bchaincat^{S_{n}}\]
Its objects are written $\mathcal{V}=\{\mathcal{V}(n)\}_{n\ge0}$.
Each collection induces an endofunctor (also denoted $\mathcal{V}$)
$\mathcal{V}:\bchaincat\to\bchaincat$\[
\mathcal{V}(X)=\bigoplus_{n\ge0}\mathcal{V}(n)\otimes_{\zs n}X^{\otimes n}\]
where $X^{\otimes n}=X\otimes\cdots\otimes X$ and $S_{n}$ acts on
$X^{\otimes n}$ by permuting factors. This endofunctor is a \emph{monad}
if the defining collection has the structure of an \emph{operad},
which means that $\mathcal{V}$ has a unit $\eta:\ring\to\mathcal{V}(1)$
and structure maps\[
\gamma_{k_{1},\dots,k_{n}}:\mathcal{V}(n)\otimes\mathcal{V}(k_{1})\otimes\cdots\otimes\mathcal{V}(k_{n})\to\mathcal{V}(k_{1}+\cdots+k_{n})\]
satisfying well-known equivariance, associativity, and unit conditions
--- see \cite{Smith-cofree}, \cite{Kriz-May}.

We will call the operad $\mathcal{V}=\{\mathcal{V}(n)\}$ $\Sigma$-\emph{cofibrant}
if $\mathcal{V}(n)$ is $\zs n$-projective for all $n\ge0$.\end{defn}
\begin{rem*}
The operads we consider here correspond to \emph{symmetric} operads
in \cite{Smith-cofree}.

The term {}``unital operad'' is used in different ways by different
authors. We use it in the sense of Kriz and May in \cite{Kriz-May},
meaning the operad has a $0$-component that acts like an arity-lowering
augmentation under compositions. Here $\mathcal{V}(0)=\ring$.

The term $\Sigma$-\emph{cofibrant} first appeared in \cite{Berger-moerdijk-axiom-operad}.
\end{rem*}
A simple example of an operad is:
\begin{example}
\label{example:frakS0}For each $n\ge0$, $C(n)=\integers S_{n}$,
with structure-map induced by\[
\gamma_{\alpha_{1},\dots,\alpha_{n}}:S_{n}\times S_{\alpha_{1}}\times\cdots\times S_{\alpha_{n}}\to S_{\alpha_{1}+\cdots+\alpha_{n}}\]
defined by regarding each of the $S_{\alpha_{i}}$ as permuting elements
within the subsequence $\{\alpha_{1}+\cdots+\alpha_{i-1}+1,\dots,\alpha_{1}+\cdots+\alpha_{i}\}$
of the sequence $\{1,\dots,\alpha_{1}+\cdots+\alpha_{n}\}$ and making
$S_{n}$ permute these $n$-blocks. This operad is denoted $\mathfrak{S}_{0}$.
In other notation, its $n^{\text{th}}$ component is the \emph{symmetric
group-ring} $\integers S_{n}$. See \cite{Smith:1994} for explicit
formulas.
\end{example}
Another important operad is:
\begin{example}
\label{exa:The-Barratt-Eccles-operad,}The Barratt-Eccles operad,
$\mathfrak{S}$, is given by $\mathfrak{S}(n)=\{C_{*}(\widetilde{K(S_{n},1})\}$
--- where $C_{*}(\widetilde{K(S_{n},1})$ is the normalized chain
complex of the universal cover of the Eilenberg-Maclane space $K(S_{n},1)$.
This is well-known (see \cite{Barratt-Eccles-operad} or \cite{Smith:1994})
to be a Hopf-operad, i.e. equipped with an operad morphism\[
\delta:\mathfrak{S}\to\mathfrak{S}\otimes\mathfrak{S}\]
and is important in topological applications. See \cite{Smith:1994}
for formulas for the structure maps.
\end{example}
For the purposes of this paper, the main example of an operad is
\begin{defn}
\label{def:coend}Given any $C\in\chaincat$, the associated \emph{coendomorphism
operad}, $\coend(C)$ is defined by\[
\coend(C)(n)=\homz(C,C^{\otimes n})\]
 Its structure map\begin{multline*}
\gamma_{\alpha_{1},\dots,\alpha_{n}}:\homz(C,C^{\otimes n})\otimes\homz(C,C^{\otimes\alpha_{1}})\otimes\cdots\otimes\homz(C,C^{\otimes\alpha_{n}})\to\\
\homz(C,C^{\otimes\alpha_{1}+\cdots+\alpha_{n}})\end{multline*}
simply composes a map in $\homz(C,C^{\otimes n})$ with maps of each
of the $n$ factors of $C$. 

This is a non-unital operad, but if $C\in\chaincat$ has an augmentation
map $\varepsilon:C\to\ring$ then we can  regard $\varepsilon$ as
the generator of $\coend(C)(0)=\ring\cdot\varepsilon\subset\homz(C,C^{\otimes0})=\homz(C,\ring)$.

Given $C\in\chaincat$ with \emph{subcomplexes} $\{D_{1},\dots,D_{k}\}$,
the \emph{relative coendomorphism operad} $\coend(C;\{D_{i}\})$ is
defined to be the sub-operad of $\coend(C)$ consisting of maps $f\in\homz(C,C^{\otimes n})$
such that $f(D_{j})\subseteq D_{j}^{\otimes n}\subseteq C^{\otimes n}$
for all $j$.
\end{defn}
We use the coendomorphism operad to define the main object of this
paper:
\begin{defn}
\label{def:coalg}A \emph{coalgebra over an operad} $\mathcal{V}$
is a chain-complex $C\in\chaincat$ with an operad morphism $\alpha:\mathcal{V}\to\coend(C)$,
called its \emph{structure map.} We will sometimes want to define
coalgebras using the \emph{adjoint structure map}\[
\alpha:C\to\prod_{n\ge0}\homzs n(\mathcal{V}(n),C^{\otimes n})\]
(in $\chaincat)$ or even the set of chain-maps\[
\alpha_{n}:C\to\homzs n(\mathcal{V}(n),C^{\otimes n})\]
for all $n\ge0$.
\end{defn}
We will sometimes want to focus on a particular class of $\mathcal{V}$-coalgebras:
the \emph{pointed, irreducible coalgebras}. We define this concept
in a way that extends the conventional definition in \cite{Sweedler:1969}:
\begin{defn}
\label{def:pointedirreducible} Given a coalgebra over a unital operad
$\mathcal{V}$ with adjoint structure-map\[
\alpha_{n}:C\to\homzs n(\mathcal{V}(n),C^{\otimes n})\]
an element $c\in C$ is called \emph{group-like} if $\alpha_{n}(c)=f_{n}(c^{\otimes n})$
for all $n>0$. Here $c^{\otimes n}\in C^{\otimes n}$ is the $n$-fold
$\ring$-tensor product, \[
f_{n}=\homz(\epsilon_{n},1):\homz(\ring,C^{\otimes n})=C^{\otimes n}\to\homzs n(\mathcal{V}(n),C^{\otimes n})\]
and $\epsilon_{n}:\mathcal{V}(n)\to\mathcal{V}(0)=\ring$ is the augmentation
(which is $n$-fold composition with $\mathcal{V}(0)$). 

A coalgebra $C$ over an operad $\mathcal{V}$ is called \emph{pointed}
if it has a \emph{unique} group-like element (denoted $1$), and \emph{pointed
irreducible} if the intersection of any two sub-coalgebras contains
this unique group-like element.\end{defn}
\begin{rem*}
Note that a group-like element generates a sub $\mathcal{V}$-coalgebra
of $C$ and must lie in dimension $0$.

Although this definition seems contrived, it arises in {}``nature'':
The chain-complex of a pointed, simply-connected reduced simplicial
set is naturally a pointed irreducible coalgebra over the Barratt-Eccles
operad, $\mathfrak{S}=\{C(K(S_{n},1))\}$ (see \cite{Smith:1994}).
In this case, the operad action encodes the chain-level effect of
Steenrod operations.\end{rem*}
\begin{prop}
Let $D$ be a pointed, irreducible coalgebra over an operad $\mathcal{V}$.
Then the augmentation map\[
\varepsilon:D\to\ring\]
is naturally split and any morphism of pointed, irreducible coalgebras
\[
f:D_{1}\to D_{2}\]
 is of the form\[
1\oplus\bar{f}:D_{1}=\ring\oplus\ker\varepsilon_{D_{1}}\to D_{2}=\ring\oplus\ker\varepsilon_{D_{2}}\]
where $\varepsilon_{i}:D_{i}\to\ring$, $i=1,2$ are the augmentations.\end{prop}
\begin{proof}
The definition (\ref{def:pointedirreducible}) of the sub-coalgebra
$\ring\cdot1\subseteq D_{i}$ is stated in an invariant way, so that
any coalgebra morphism must preserve it. Any morphism must also preserve
augmentations because the augmentation is the $0^{\mathrm{th}}$-order
structure-map. Consequently, $f$ must map $\ker\varepsilon_{D_{1}}$to
$\ker\varepsilon_{D_{2}}$. The conclusion follows.\end{proof}
\begin{defn}
\label{def:pointedirredcat} We denote the \emph{category} of nearly
free coalgebras over $\mathcal{V}$ by $\coalgcat$. The terminal
object in this category is $0$, the null coalgebra.

The category of nearly free \emph{pointed irreducible coalgebras}
over $\mathcal{V}$ is denoted $\ircoalgcat$ --- this is only defined
if $\mathcal{V}$ is unital. \emph{Its} terminal object is the coalgebra
whose underlying chain complex is $\ring$ concentrated in dimension
$0$ with coproduct that sends $1\in\ring$ to $1^{\otimes n}\in\ring^{\otimes n}$.

It is not hard to see that these terminal objects are also the \emph{initial}
objects of their respective categories.
\end{defn}
We also need:
\begin{defn}
\label{def:forgetful}If $A\in\mathrf C=\ircoalgcat$ or $\coalgcat$,
then $\forgetful A$ denotes the underlying chain-complex in $\chaincat$
of\[
\ker A\to\bullet\]
where $\bullet$ denotes the terminal object in $\mathrf{C}$ ---
see definition~ \ref{def:pointedirredcat}. We will call $\forgetful{\ast}$
the \emph{forgetful functor} from $\mathrf{C}$ to $\chaincat$.
\end{defn}
We can also define the analogue of an ideal:
\begin{defn}
\label{def:coideal}Let $C$ be a coalgebra over the operad $\mathcal{U}$
with adjoint structure map \[
\alpha:C\to\prod_{n\ge0}\homzs n(\mathcal{U}(n),C^{\otimes n})\]
and let $D\subseteq\forgetful C$ be a sub-chain complex that is a
direct summand. Then $D$ will be called a \emph{coideal} of $C$
if the composite \[
\alpha|D:D\to\prod_{n\ge0}\homzs n(\mathcal{U}(n),C^{\otimes n})\xrightarrow{\homz(1_{\mathcal{U}},p^{\otimes})}\prod_{n\ge0}\homzs n(\mathcal{U}(n),(C/D)^{\otimes n})\]
 \emph{vanishes}, where $p:C\to C/D$ is the projection to the quotient
(in $\chaincat$). \end{defn}
\begin{rem*}
Note that it is easier for a sub-chain-complex to be a coideal of
a coalgebra than to be an ideal of an algebra. For instance, all sub-coalgebras
of a coalgebra are also coideals. Consequently it is easy to form
quotients of coalgebras and hard to form sub-coalgebras. This is dual
to what occurs for algebras.
\end{rem*}
We will use the concept of cofree coalgebra cogenerated by a chain
complex: 
\begin{defn}
\label{def:cofreecoalgebra}Let $C\in\chaincat$ and let $\mathcal{V}$
be an operad. Then a $\mathcal{V}$-coalgebra $G$ will be called
\emph{the cofree coalgebra} \emph{cogenerated by} $C$ if
\begin{enumerate}
\item there exists a morphism of DG-modules $\varepsilon:G\to C$
\item given any $\mathcal{V}$-coalgebra $D$ and any morphism of DG-modules$f:D\to C$,
there exists a \emph{unique} morphism of $\mathcal{V}$-coalgebras,
$\hat{f}:D\to G$, that makes the diagram \[
\xymatrix{{D}\ar[r]^{\hat{f}}\ar[rd]_{f} & {G}\ar[d]^{\varepsilon}\\
 & {C}}
\]
 commute.
\end{enumerate}
\end{defn}
This universal property of cofree coalgebras implies that they are
unique up to  isomorphism if they exist. The paper \cite{Smith-cofree}
gives a constructive proof of their existence in great generality
(under the unnecessary assumption that chain-complexes are $\ring$-free).
In particular, this paper defines cofree coalgebras $L_{\mathcal{V}}C$
and pointed irreducible cofree coalgebras $P_{\mathcal{V}}C$ cogenerated
by a chain-complex $C$. There are several ways to define them:
\begin{enumerate}
\item $L_{\mathcal{V}}C$ is essentially the largest submodule of \[
C\oplus\prod_{k=1}^{\infty}\homzs k(\mathcal{V}(k),C^{\otimes k})\]
on which the coproduct defined by the dual of the composition-operations
of $\mathcal{V}$ is well-defined.
\item If $C$ is a coalgebra over $\mathcal{V}$, its image under the structure
map\[
C\to C\oplus\prod_{k=1}^{\infty}\homzs k(\mathcal{V}(k),C^{\otimes k})\]
turns out to be a sub-coalgebra of the target --- with a coalgebra
structure that vanishes on the left summand ($C$) and is the dual
of the structure-map of $\mathcal{V}$ on the right. We may define
$L_{\mathcal{V}}C$ to be the sum of all coalgebras in $C\oplus\prod_{k=1}^{\infty}\homzs k(\mathcal{V}(k),C^{\otimes k})$
formed in this way. The classifying map of a coalgebra \[
C\to L_{\mathcal{V}}C\]
 is just the structure map of the coalgebra structure.
\end{enumerate}

\section{Model categories\label{sec:Model-categories}}

We recall the concept of a model structure on a category $\mathcal{G}$.
This involves defining specialized classes of morphisms called \emph{cofibrations,}
\emph{fibrations,} and \emph{weak equivalences} (see \cite{Quillen:1967}
and \cite{Goerss-Jardine}). The category and these classes of morphisms
must satisfy the conditions:
\begin{description}
\item [{CM~1}] $\mathcal{G}$ is closed under all finite limits and colimits
\item [{CM~2}] Suppose the following diagram commutes in $\mathcal{G}$:
\[
\xymatrix{{X}\ar[rr]^{g}\ar[dr]_{h} & {} & {Y}\\
{} & {Z}\ar[ur]_{f} & {}}
\]
If any two of $f,g,h$ are weak equivalences, so is the third.
\item [{CM~3}] These classes of morphisms are closed under formation of
\emph{retracts:} Given a commutative diagram \[
\xymatrix{{A}\ar[r]\ar[d]_{f} & {B}\ar[r]\ar[d]_{g} & {A}\ar[d]^{f}\\
{C}\ar[r] & {D}\ar[r] & {C}}
\]
 whose horizontal composites are the identity map, if $g$ is a weak
equivalence, fibration, or cofibration, then so is $f$. 
\item [{CM~4}] Given a commutative solid arrow diagram \[
\xymatrix{{U}\ar[r]\ar[d]_{i} & {A}\ar[d]^{p}\\
{W}\ar[r]\ar@{.>}[ur] & {B}}
\]
 where $i$ is a cofibration and $p$ is a fibration, the dotted arrow
exists whenever $i$ or $p$ are trivial. 
\item [{CM~5}] Any morphism $f:X\to Y$ in $\mathcal{G}$ may be factored:

\begin{enumerate}
\item $f=p\circ i$, where $p$ is a fibration and $i$ is a trivial cofibration
\item $f=q\circ j$, where $q$ is a trivial fibration and $j$ is a cofibration
\end{enumerate}
We also assume that these factorizations are \emph{functorial} ---
see \cite{Homotopy-theories-and-model-cats}.

\end{description}
\begin{defn}
An object, $X$, for which the map $\bullet\to X$ is a cofibration,
is called \emph{cofibrant}. An object, $Y$, for which the map $Y\to\bullet$
is a fibration, is called \emph{fibrant}.

The properties of a model category immediately imply that:\end{defn}
\begin{lem}
\label{lem:trivfibrationretraction}Let $f:A\to B$ be a trivial fibration
in which $B$ is cofibrant. Then $f$ is a retraction of $A$ onto
$B$, i.e., there exists a morphism $g:B\to A$ such that $f\circ g=1:B\to B$.\end{lem}
\begin{proof}
Consider the diagram \[
\xymatrix{{\bullet}\ar[r]\ar[d] & {A}\ar[d]^{f}\\
{B}\ar@{=}[r] & {B}}
\]
 Property CM~4 implies that we can complete this to a diagram \[
\xymatrix{{\bullet}\ar[r]\ar[d] & {A}\ar[d]^{f}\\
{B}\ar@{.>}[ru]^{g}\ar@{=}[r] & {B}}
\]

\end{proof}

\subsection{A model-category of chain-complexes\label{ex:absolutemodel}}

Let $\chaincat$ denote the category of unbounded chain-complexes
over the ring $\ring$. The \emph{absolute model structure} of Christensen
and Hovey in \cite{Christensen.Hovey:2002}, and Cole in \cite{Cole.M-homotopy}
is defined via:
\begin{enumerate}
\item \emph{Weak equivalences} are chain-homotopy equivalences: two chain-complexes
$C$ and $D$ are weakly equivalent if there exist chain-maps: $f:C\to D$
and $g:D\to C$ and chain-homotopies $\varphi_{1}:C\to C$ and $\varphi_{2}:D\to D$
such that $d\varphi_{1}=g\circ f-1$ and $d\varphi_{2}=f\circ g-1$. 
\item \emph{Fibrations} are surjections of chain-complexes that are split
(as maps of graded $\ring$-modules).
\item \emph{Cofibrations} are injections of chain-complexes that are split
(as maps of graded $\ring$-modules).\end{enumerate}
\begin{rem*}
\emph{All} chain complexes are fibrant and cofibrant in this model.

In this model structure, a quasi-isomorphism may \emph{fail} to be
a weak equivalence. It is well-known not to be cofibrantly generated
(see \cite{Christensen.Hovey:2002}).

Since all chain-complexes are cofibrant, lemma~\ref{lem:trivfibrationretraction}
implies that all trivial fibrations are retractions --- i.e., they
are split as chain-maps. We will need the following relative version
of lemma~\ref{lem:trivfibrationretraction} in the sequel:\end{rem*}
\begin{lem}
\label{lem:reltrivfibrationsplit}Let \[
\xymatrix{{A}\ar[r]^{u}\ar[d]_{f} & {B}\ar[d]^{g}\\
{C}\ar[r]_{v} & {D}}
\]
 be a commutative diagram in $\chaincat$ such that 
\begin{itemize}
\item $f$ and $g$ are trivial fibrations
\item $v$ is a cofibration
\end{itemize}
Then there exists maps $\ell:C\to A$ and $m:D\to B$ such that 
\begin{enumerate}
\item $f\circ\ell=1:C\to C$, $g\circ m=1:D\to D$ 
\item the diagram\[
\xymatrix{{A}\ar[r]^{u} & {B}\\
{C}\ar[r]_{v}\ar[u]^{\ell} & {D}\ar[u]_{m}}
\]
commutes.
\end{enumerate}
If $u$ is injective and split, there exist homotopies $\varphi_{1}:A\otimes I\to A$
from $1$ to $\ell\circ f$ and $\varphi_{2}:B\otimes I\to B$ from
$1$ to $m\circ g$, respectively, such that the diagram\[
\xymatrix{{A\otimes I}\ar[d]_{\varphi_{1}}\ar[r]^{u\otimes1} & {B\otimes I}\ar[d]^{\varphi_{2}}\\
{A}\ar[r]_{u} & {B}}
\]
commutes.\end{lem}
\begin{proof}
We construct $\ell$ exactly as in lemma~\ref{lem:trivfibrationretraction}
and use it to create a commutative diagram \[
\xymatrix{{C}\ar[r]^{u\circ\ell}\ar[d]_{v} & {B}\ar[d]^{g}\\
{D}\ar@{=}[r] & {D}}
\]
Property CM~4 implies that we can complete this to a commutative
diagram\[
\xymatrix{{C}\ar[r]^{u\circ\ell}\ar[d]_{v} & {B}\ar[d]^{g}\\
{D}\ar@{=}[r]\ar@{.>}[ru]^{m} & {D}}
\]
Let $\varphi_{1}':A\otimes I\to A$ and $\varphi_{2}':B\otimes I\to B$
be any homotopies from the identity map to $\ell\circ f$ and $m\circ g$,
respectively. If there exists a chain-map $p:B\to A$ $p\circ u=1:A\to A$
then $\varphi_{1}=\varphi_{1}'$ and $\varphi_{2}=u\circ\varphi_{1}'\circ(p\otimes1)+\varphi_{2}'\circ(1-(u\circ p)\otimes1)$
have the required properties.
\end{proof}

\section{Model-categories of coalgebras\label{sec:The-general-case}}

\subsection{Description of the model-structure\label{sub:Description-of-the}}

We will base our model-structure on that of the underlying chain-complexes
in $\chaincat$. Definition~\ref{def:cofibrationweakequiv} and definition
\ref{def:fibration} describe how we define cofibrations, fibrations,
and weak equivalences.

We must allow non-$\ring$-free chain-complexes (when $\ring=\integers$)
because the underlying chain complexes of the cofree coalgebras $\pcoalg V{\ast}$
and $L_{\mathcal{V}}(\ast)$ are not known to be $\ring$-free. They
certainly are if $\ring$ is a field, but if $\ring=\integers$ their
underlying abelian groups are subgroups of the Baer-Specker group,
$\integers^{\aleph_{0}}$, which is $\integers$-torsion free but
well-known \emph{not} to be a free abelian group (see \cite{Baer-Specker-nonfree},
\cite{Baer-Specker-bounded} or the survey \cite{Baer-Specker-survey}). 
\begin{prop}
The forgetful functor (defined in definition~\ref{def:forgetful})
and cofree coalgebra functors define adjoint pairs\begin{eqnarray*}
\pcoalg V{\ast}:\chaincat & \leftrightarrows & \ircoalgcat:\forgetful{\ast}\\
L_{\mathcal{V}}(\ast):\chaincat & \leftrightarrows & \coalgcat:\forgetful{\ast}\end{eqnarray*}
\end{prop}
\begin{rem*}
The adjointness of the functors follows from the universal property
of cofree coalgebras --- see~\cite{Smith-cofree}. \end{rem*}
\begin{condition}
\emph{\label{cond:mainoperadassumption}Throughout the rest of this
paper,} we assume that $\mathcal{V}$ is an operad equipped with a
morphism of operads\[
\delta:\mathcal{V}\to\mathcal{V}\otimes\mathfrak{S}\]
--- where $\mathfrak{S}$ is the Barratt-Eccles operad (see example~\ref{exa:The-Barratt-Eccles-operad,}
--- that makes the diagram \[
\xymatrix{{\mathcal{V}}\ar[r]^{\delta}\ar@{=}[rd] & {\mathcal{V}\otimes\mathfrak{S}}\ar[d]\\
{} & {\mathcal{V}}}
\]
 commute. Here, the operad structure on $\mathcal{V}\otimes\mathfrak{S}$
is just the tensor product of the operad structures of $\mathcal{V}$
and $\mathfrak{S}$, and the vertical map is projection:\[
\mathcal{V}\otimes\mathfrak{S}\xrightarrow{1\otimes\hat{\epsilon}}\mathcal{V}\otimes T=\mathcal{V}\]
where $T$ is the operad that is $\ring$ in all arities and $\hat{\epsilon}:\mathfrak{S}\to T$
is defined by the augmentations:\[
\epsilon_{n}:\rs n\to\ring\]

In addition, we assume that, for each $n\ge0$, $\{\mathcal{V}(n)\}$
is an $\zs n$-projective chain-complex of finite type. 

We also assume that the arity-$1$ component of $\mathcal{V}$ is
equal to $\ring$, generated by the unit.\end{condition}
\begin{rem}
Free and cofibrant operads (with each component of finite type) satisfy
this condition. The condition that the chain-complexes are projective
corresponds to the Berger and Moerdijk's condition of $\Sigma$-cofibrancy
in \cite{Berger-moerdijk-axiom-operad}.
\end{rem}
Now we define our model structure on the categories $\ircoalgcat$
and $\coalgcat$. 
\begin{defn}
\label{def:cofibrationweakequiv}A morphism $f:A\to B$ in $\mathrf{C}=\coalgcat$
or $\ircoalgcat$ will be called
\begin{enumerate}
\item a \emph{weak equivalence} if $\forgetful f:\forgetful A\to\forgetful B$
is a chain-homotopy equivalence in $\chaincat$. An object $A$ will
be called \emph{contractible} if the augmentation map\[
A\to\bullet\]
is a weak equivalence, where $\bullet$ denotes the terminal object
in $\mathrf{C}$ --- see definition~\ref{def:pointedirredcat}. 
\item a \emph{cofibration} if $\forgetful f$ is a cofibration in $\chaincat$.
\item a \emph{trivial cofibration} if it is a weak equivalence and a cofibration.
\end{enumerate}
\end{defn}
\begin{rem*}
A morphism is a cofibration if it is a degreewise split monomorphism
of chain-complexes. Note that all objects of $\mathrf{C}$ are cofibrant.

Our definition makes $f:A\to B$ a weak equivalence if and only if
$\forgetful f:\forgetful A\to\forgetful B$ is a weak equivalence
in $\chaincat$. \end{rem*}
\begin{defn}
\label{def:fibration}A morphism $f:A\to B$ in $\coalgcat$ or $\ircoalgcat$
will be called 
\begin{enumerate}
\item a \emph{fibration} if the dotted arrow exists in every diagram of
the form \[
\xymatrix{{U}\ar[r]\ar[d]_{i} & {A}\ar[d]^{f}\\
{W}\ar[r]\ar@{.>}[ur] & {B}}
\]
 in which $i:U\to W$ is a trivial cofibration.
\item a \emph{trivial fibration} if it is a fibration and a weak equivalence.
\end{enumerate}
\end{defn}
Definition~\ref{def:cofibrationweakequiv} explicitly described cofibrations
and definition~\ref{def:fibration} defined fibrations in terms of
them. We will verify the axioms for a model category (part of CM~4
and CM~5) and characterize fibrations.

We will occasionally need a stronger form of equivalence:
\begin{defn}
\label{def:strictequivalence}Let $f,g:A\to B$ be a pair of morphisms
in $\coalgcat$ or $\ircoalgcat$. A \emph{strict homotopy} between
them is a coalgebra-morphism (where $A\otimes I$ has the coalgebra
structure defined in condition~\ref{cond:mainoperadassumption})\[
F:A\otimes I\to B\]
such that $F|A\otimes p_{0}=f:A\otimes p_{0}\to B$ and $F|A\otimes p_{1}=g:A\otimes p_{01}\to B$.
A \emph{strict equivalence} between two coalgebras $A$ and $B$ is
a pair of coalgebra-morphisms\begin{eqnarray*}
f:A & \to & B\\
g:B & \to & A\end{eqnarray*}
and strict homotopies from $f\circ g$ to the identity of $B$ and
from $g\circ f$ to the identity map of $A$.\end{defn}
\begin{rem*}
Strict equivalence is a direct translation of the definition of weak
equivalence in $\chaincat$ into the realm of coalgebras. Strict equivalences
are weak equivalences but the converse is not true.

The reader may wonder why we didn't use strict equivalence in place
of what is defined in definition~\ref{def:cofibrationweakequiv}.
It turns out that in we are only able to prove CM~5 with the weaker
notion of equivalence used here. 
\end{rem*}
In a few simple cases, describing fibrations is easy: 
\begin{prop}
\label{pro:surjectivecofreefibration}Let \[
f:A\to B\]
be a fibration in $\chaincat$. Then the induced morphisms\begin{eqnarray*}
P_{\mathcal{V}}f:P_{\mathcal{V}}A & \to & P_{\mathcal{V}}B\\
L_{\mathcal{V}}f:L_{\mathcal{V}}A & \to & L_{\mathcal{V}}B\end{eqnarray*}
are fibrations in $\ircoalgcat$ and $\coalgcat$, respectively.\end{prop}
\begin{proof}
Consider the diagram \[
\xymatrix{{U}\ar[r]\ar[d] & {P_{\mathcal{V}}A}\ar[d]^{P_{\mathcal{V}}f}\\
{V}\ar[r]\ar@{.>}[ur] & {P_{\mathcal{V}}B}}
\]
 where $U\to V$ is a trivial cofibration --- i.e., $\forgetful U\to\forgetful V$
is a trivial cofibration of \emph{chain-complexes.} Then the dotted
map exists by the the defining property of cofree coalgebras and by
the existence of the lifting map in the diagram \[
\xymatrix{{\forgetful U}\ar[r]\ar[d] & {A}\ar[d]^{f}\\
{\forgetful V}\ar[r]\ar@{.>}[ur] & {B}}
\]
 of chain-complexes.\end{proof}
\begin{cor}
\label{cor:cofreefibrant}All cofree coalgebras are fibrant.\end{cor}
\begin{prop}
\label{prop:lefthomotopy}Let $C$ and $D$ be objects of $\chaincat$
and let \[
f_{1},f_{2}:C\to D\]
 be chain-homotopic morphisms via a chain-homotopy \begin{equation}
F:C\otimes I\to D\label{eq:chainhomotop}\end{equation}
Then the induced maps\begin{eqnarray*}
P_{\mathcal{V}}f_{i}:P_{\mathcal{V}}C & \to & P_{\mathcal{V}}D\\
L_{\mathcal{V}}f_{i}:L_{\mathcal{V}}C & \to & L_{\mathcal{V}}D\end{eqnarray*}
$i=1,2$, are left-homotopic in $\ircoalgcat$ and $\coalgcat$, respectively
via a strict chain homotopy\[
F':P_{\mathcal{V}}f_{i}:(P_{\mathcal{V}}C)\otimes I\to P_{\mathcal{V}}D\]
If we equip $C\otimes I$ with a coalgebra structure using condition~\ref{cond:mainoperadassumption}
and if $F$ in \ref{eq:chainhomotop} is strict then the diagram \[
\xymatrix{{C\otimes I}\ar[r]^{F}\ar[d]_{\alpha_{C}\otimes1} & {D}\ar[d]^{\alpha_{D}}\\
{P_{\mathcal{V}}(C)\otimes I}\ar[r]_{\quad F'} & {P_{\mathcal{V}}D}}
\]
 commutes in the pointed irreducible case and the diagram \[
\xymatrix{{C\otimes I}\ar[r]^{F}\ar[d]_{\alpha_{C}\otimes1} & {D}\ar[d]^{\alpha_{D}}\\
{L_{\mathcal{V}}(C)\otimes I}\ar[r]_{\quad F'} & {L_{\mathcal{V}}D}}
\]
commutes in the general case. Here $\alpha_{C}$ and $\alpha_{D}$
are classifying maps of coalgebra structures.\end{prop}
\begin{rem*}
In other words, the cofree coalgebra functors map homotopies and weak
equivalences in $\chaincat$ to strict homotopies and strict equivalences,
respectively, in $\ircoalgcat$ and $\coalgcat$.

If the homotopy in $\chaincat$ was the result of applying the forgetful
functor to a strict homotopy, then the generated strict homotopy is
compatible with it.\end{rem*}
\begin{proof}
We will prove this in the pointed irreducible case. The general case
follows by a similar argument. The chain-homotopy between the $f_{i}$
induces\[
P_{\mathcal{V}}F:P_{\mathcal{V}}(C\otimes I)\to P_{\mathcal{V}}D\]
Now we construct the map \[
H:(P_{\mathfrak{\mathcal{V}}}C)\otimes I\to P_{\mathfrak{\mathcal{V}}}(C\otimes I)\]
using the universal property of a cofree coalgebra and the fact that
the coalgebra structure of $(P_{\mathcal{V}}C)\otimes I$ extends
that of $P_{\mathfrak{\mathcal{V}}}C$ on both ends by condition~\ref{cond:mainoperadassumption}.
Clearly\[
P_{\mathcal{V}}F\circ H:(P_{\mathcal{V}}C)\otimes I\to P_{\mathcal{V}}D\]
is the required left-homotopy.

If we define a coalgebra structure on $C\otimes I$ using condition~\ref{cond:mainoperadassumption},
we get diagram \[
\xymatrix{{C\otimes I}\ar@{=}[r]\ar[d]_{\alpha_{C}\otimes1} & {C\otimes I}\ar[r]^{\quad F}\ar[d]^{\alpha_{C\otimes I}} & {D}\ar[d]^{\alpha_{D}}\\
{P_{\mathcal{V}}(C)\otimes I}\ar[d]_{\epsilon_{C}\otimes1}\ar[r]^{H} & {P_{\mathcal{V}}(C\otimes I)}\ar[r]_{\quad P_{\mathcal{V}}F}\ar[d]^{\epsilon_{C\otimes I}} & {P_{\mathcal{V}}D}\\
{C\otimes I}\ar@{=}[r] & {C\otimes I} & {}}
\]
where $\alpha_{C\otimes I}$ is the classifying map for the coalgebra
structure on $C\otimes I$. 

We claim that this diagram commutes. The fact that $F$ is a coalgebra
morphism implies that the upper right square commutes. The large square
on the left (bordered by $C\otimes I$ on all four corners) commutes
by the property of co-generating maps (of cofree coalgebras) and classifying
maps. The two smaller squares on the left (i.e., the large square
with the map $H$ added to it) commute by the universal properties
of cofree coalgebras (which imply that induced maps to cofree coalgebras
are uniquely determined by their composites with co-generating maps).
The diagram in the statement of the result is just the outer upper
square of this diagram, so we have proved the claim.
\end{proof}
This result implies a homotopy invariance property of the \emph{categorical
product,} $A_{0}\boxtimes A_{1}$, defined explicitly in definition~\ref{def:catprod}
of appendix~\ref{sec:Basic-constructions}. 
\begin{lem}
\label{lem:pullbackfibrationfibration}Let $g:B\to C$ be a fibration
in $\ircoalgcat$ and let $f:A\to C$ be a morphism in $\ircoalgcat$.
Then the projection\[
A\boxtimes^{C}B\to A\]
is a fibration.\end{lem}
\begin{rem*}
The notation $A\boxtimes^{C}B$ denotes a fibered product --- see
definition~\ref{def:catfiberedproduct} in appendix~\ref{sub:Categorical-products}
for the precise definition. In other words, pullbacks of fibrations
are fibrations.\end{rem*}
\begin{proof}
Consider the diagram \begin{equation}
\xymatrix{{U}\ar[r]^{u\quad\quad}\ar[d]_{i} & {A\boxtimes^{C}B}\ar[d]^{p_{A}}\\
{V}\ar[r]_{v}\ar@{.>}[ur] & {A}}
\label{dia:text1}\end{equation}
 where $U\to V$ is a trivial cofibration. The defining property of
a categorical product implies that any map to $A\boxtimes^{C}B\subseteq A\boxtimes B$
is determined by its composites with the projections \begin{eqnarray*}
p_{A}:A\boxtimes B & \to & A\\
p_{B}:A\boxtimes B & \to & B\end{eqnarray*}
Consider the composite $p_{B}\circ u:U\to B$. The commutativity of
the solid arrows in diagram~\ref{dia:text1} implies that the diagram\[
\xymatrix{{U}\ar[rr]^{p_{B}\circ u}\ar[dr]^{p_{A}\circ u}\ar[d]_{i} & {} & {B}\ar[d]^{g}\\
{V}\ar[r]_{v} & {A}\ar[r]_{f} & {C}}
\]
 commutes and this implies that the solid arrows in the diagram \begin{equation}
\xymatrix{{U}\ar[r]^{p_{B}\circ u}\ar[d]_{i} & {B}\ar[d]^{g}\\
{V}\ar[r]_{f\circ v}\ar@{.>}[ur] & {C}}
\label{dia:text2}\end{equation}
 commute. The fact that $g:B\to C$ is a \emph{fibration} implies
that the dotted arrow exists in diagram~\ref{dia:text2}, which implies
the existence of a map $V\to A\boxtimes B$ whose composites with
$f$ and $g$ agree. This defines a map $V\to A\boxtimes^{C}B$ that
makes \emph{all} of diagram~\ref{dia:text1} commute. The conclusion
follows.
\end{proof}

\subsection{Proof of CM~1 through CM~3\label{sub:Proof-of-CM1}}

CM~1 asserts that our categories have all finite limits and colimits. 

The results of appendix~\ref{sec:Basic-constructions} prove that
all \emph{countable} limits and colimits exist --- see theorem~\ref{th:cofreelimits}
and theorem~\ref{th:coalgdirectlimit}.

CM~2 follows from the fact that we define weak equivalence the same
was it is defined in $\chaincat$ --- so the model structure on $\chaincat$
implies that this condition is satisfies on our categories of coalgebras.
A similar argument verifies condition CM~3.

\subsection{Proof of CM~5\label{sub:Proof-of-CM5}}

We begin with:
\begin{cor}
Let $A\in\ircoalgcat$ be fibrant and let $B\in\ircoalgcat$. Then
the projection\[
A\boxtimes B\to B\]
is a fibration.
\end{cor}
This allows us to verify CM 5, statement 2:
\begin{cor}
\label{cor:trivcofibfibfactor}Let $f:A\to B$ be a morphism in $\mathrf{C}=\ircoalgcat$
or $\coalgcat$, and let\[
Z=\left\{ \begin{array}{ll}
P_{\mathcal{V}}\cone{\forgetful A}\boxtimes B & \text{when\,\,\,}\mathrf{C}=\ircoalgcat\\
L_{\mathcal{V}}\cone{\forgetful A}\boxtimes B & \text{when\,\,\,}\mathrf{C}=\coalgcat\end{array}\right.\]
Then $f$ factors as\[
A\to Z\to B\]
where 
\begin{enumerate}
\item $\cone{\forgetful A}$ is the cone on $\forgetful A$ (see definition~\ref{def:unitinterval})
with the canonical inclusion $i:\forgetful A\to\cone{\forgetful A}$ 
\item the morphism $i\boxtimes f:A\to Z$ is a cofibration
\item the morphism $Z\to B$ is projection to the second factor and is a
fibration.
\end{enumerate}
Consequently, $f$ factors as a cofibration followed by a trivial
fibration.

\end{cor}
\begin{proof}
We focus on the pointed irreducible case. The general case follows
by essentially the same argument. The existence of the (injective)
morphism $A\to P_{\mathcal{V}}\cone{\forgetful A}\boxtimes B$ follows
from the definition of $\boxtimes$. We claim that its image is a
direct summand of $P_{\mathcal{V}}\cone{\forgetful A}\boxtimes B$
as a graded $\ring$-module (which implies that $i\boxtimes f$ is
a cofibration). We clearly get a projection\[
P_{\mathcal{V}}\cone{\forgetful A}\boxtimes B\to P_{\mathcal{V}}\cone{\forgetful A}\]
 and the composite of this with the co-generating map $\forgetful{P_{\mathcal{V}}\cone{\forgetful A}}\to\cone{\forgetful A}$
gives rise a a morphism of chain-complexes\begin{equation}
\forgetful{P_{\mathcal{V}}\cone{\forgetful A}\boxtimes B}\to\cone{\forgetful A}\label{eq:splittingchaincomplex}\end{equation}
Now note the existence of a splitting map\[
\cone{\forgetful A}\to\forgetful A\]
of \emph{graded $\ring$-modules} (not coalgebras or \emph{even} chain-complexes).
Combined with the map in equation~\ref{eq:splittingchaincomplex},
we conclude that $A\to P_{\mathcal{V}}\cone{\forgetful A}\boxtimes B$
is a cofibration.

There is a weak equivalence $c:\cone{\forgetful A}\to\bullet$ in
$\chaincat$, and \ref{prop:lefthomotopy} implies that it induces
a strict equivalence $ $$P_{\mathcal{V}}c:P_{\mathcal{V}}\cone{\forgetful A}\to\bullet$.
Proposition~\ref{pro:stricthomotopyproducts} implies that \[
c\boxtimes1:P_{\mathcal{V}}\cone{\forgetful A}\boxtimes B\to\bullet\boxtimes B=B\]
 is a strict equivalence.
\end{proof}
The first part of CM~5 will be considerably more difficult to prove.
\begin{defn}
\label{def:pro-categories}Let $\iircoalgcat$ and $\icoalgcat$ be
the categories of inverse systems of objects of $\ircoalgcat$ and
$\coalgcat$, respectively and let $\dircoalgcat$ and $\dcoalgcat$
be corresponding categories of direct systems. Morphisms are defined
in the obvious way.
\end{defn}
Now we define the \emph{rug-resolution} of a cofibration:
\begin{defn}
\label{def:rug-resolution}Let $\mathcal{V}=\{\mathcal{V}(n)\}$ be
a $\Sigma$-cofibrant (see definition~\ref{def:operad}) operad such
that $\mathcal{V}(n)$ is of finite type for all $n\ge0$. If $f:C\to D$
is a cofibration in $\ircoalgcat$ or $\coalgcat$, define\begin{eqnarray*}
G_{0} & = & D\\
f_{0}=f:C & \to & G_{0}\\
G_{n+1} & = & G_{n}\boxtimes^{L_{\mathcal{V}}\forgetful{H_{n}}}L_{\mathcal{V}}\bar{H}_{n}\\
p_{n+1}:G_{n+1} & \to & G_{n}\end{eqnarray*}
 for all $n$, where
\begin{enumerate}
\item $\epsilon:C\to\bullet$ is the unique morphism.
\item $H_{n}$ is the cofiber of $f_{n}$ in the push-out\[
\xymatrix{C\ar[r]^{\epsilon}\ar[d]_{f_{n}} & {\bullet}\ar@{.>}[d]\\
{G_{n}}\ar@{.>}[r] & {H_{n}}}
\]

\item $G_{n}\to L_{\mathcal{V}}\forgetful{H_{n}}$ is the composite of the
classifying map \[
G_{n}\to L_{\mathcal{V}}\forgetful{G_{n}}\]
with the map\[
L_{\mathcal{V}}\forgetful{G_{n}}\to L_{\mathcal{V}}\forgetful{H_{n}}\]

\item $\bar{H_{n}}=\Sigma^{-1}\cone{\forgetful{H_{n}}}$ --- where $\Sigma^{-1}$
denotes desuspension (in $\chaincat$). It is contractible and comes
with a canonical $\chaincat$-fibration\begin{equation}
v_{n}:\bar{H_{n}}\to\forgetful{H_{n}}\label{eq:vndef}\end{equation}
inducing the fibration\begin{equation}
L_{\mathcal{V}}v_{n}:L_{\mathcal{V}}\bar{H}_{n}\to L_{\mathcal{V}}\forgetful{H_{n}}\label{eq:lvn}\end{equation}

\item $p_{n+1}:G_{n+1}=G_{n}\boxtimes^{L_{\mathcal{V}}\forgetful{H_{n}}}L_{\mathcal{V}}\bar{H}_{n}\to G_{n}$
is projection to the first factor,
\item The map $f_{n+1}:C\to G_{n}\boxtimes^{L_{\mathcal{V}}\forgetful{H_{n}}}L_{\mathcal{V}}\bar{H}_{n}$
is the \emph{unique} morphism that makes the diagram \[
\xymatrix{{} & {G_{n}\boxtimes L_{\mathcal{V}}\bar{H}_{n}}\ar[ld]\ar[rd] & {}\\
{G_{n}} & {} & {L_{\mathcal{V}}\bar{H}_{n}}\\
{} & {C}\ar[lu]^{f_{n}}\ar[ru]_{\epsilon}\ar@{.>}[uu]_{f_{n+1}} & {}}
\]
commute, where the downwards maps are projections to factors. The
map $\epsilon:C\to L_{\mathcal{V}}\bar{H}_{n}$ is

\begin{enumerate}
\item the map to the \emph{basepoint} if the category is $\ircoalgcat$
(and $L_{\mathcal{V}}\bar{H}_{n}$ is replaced by $P_{\mathcal{V}}\bar{H}_{n}$),
\item the \emph{zero}-map if the category is $\coalgcat$.
\end{enumerate}
The commutativity of the diagram \[
\xymatrix{{} & L_{\mathcal{V}}H_{n} & {}\\
{G_{n}}\ar[ru] & {} & {L_{\mathcal{V}}\bar{H}_{n}}\ar[lu]\\
{} & {C}\ar[lu]^{f_{n}}\ar[ru]_{\epsilon} & {}}
\]
implies that the image of $f_{n+1}$ actually lies in the fibered
product, $G_{n}\boxtimes^{L_{\mathcal{V}}\forgetful{H_{n}}}L_{\mathcal{V}}\bar{H}_{n}$.

\end{enumerate}
The \emph{rug-resolution} of $f:C\to D$ is the map of inverse systems
$\{f_{i}\}:\underline{\{C\}}\to\{G_{i}\}\to D$, where $\underline{\{C\}}$
denotes the constant inverse system.\end{defn}
\begin{rem*}
Very roughly speaking, this produces something like a {}``Postnikov
resolution'' for $f:C\to D$. Whereas a Postnikov resolution's stages
{}``push the trash upstairs,'' this one's {}``push the trash horizontally''
or {}``under the rug'' --- something feasible because one has an
infinite supply of rugs. \end{rem*}
\begin{prop}
\label{pro:hnnummhomotopic}Following all of the definitions of \ref{def:rug-resolution}
above, the diagrams \[
\xymatrix{{} & {G_{n+1}}\ar[d]^{p_{n+1}}\\
{C}\ar[ru]^{f_{n+1}}\ar[r]_{f_{n}} & {G_{n}}}
\]
commute and induce maps $H_{n+1}\xrightarrow{p'_{n+1}}H_{n}$ that
fit into commutative diagrams of chain-complexes\[
\xymatrix{{} & {\bar{H}_{n}}\ar[d]^{e_{n}}\\
{H_{n+1}}\ar[ru]^{u_{n}}\ar[r]_{p'_{n+1}} & {H_{n}}}
\]
It follows that the maps $p'_{n+1}$ are nullhomotopic for all $n$.\end{prop}
\begin{proof}
Commutativity is clear from the definition of $f_{n+1}$ in terms
of $f_{n}$ above. 

To see that the induced maps are nullhomotopic, consider the diagram
\[
\xymatrix{{G_{n}\boxtimes^{L_{\mathcal{V}}\forgetful{H_{n}}}L_{\mathcal{V}}\bar{H}_{n}}\ar[rr]\ar[d]_{p_{n+1}} & {} & {L_{\mathcal{V}}\bar{H}_{n}}\ar[r]^{\varepsilon}\ar[d]_{L_{\mathcal{V}}v_{n}} & {\bar{H}_{n}}\ar[d]^{v_{n}}\\
{G_{n}}\ar[r] & {H_{n}}\ar[r]_{\alpha} & {L_{\mathcal{V}}H_{n}}\ar[r]_{\varepsilon} & {H_{n}}}
\]
where $v_{n}$ is defined in equation~\ref{eq:vndef}, both $\varepsilon$-maps
are cogenerating maps --- see definition~\ref{def:cofreecoalgebra}
--- and $\alpha:H_{n}\to L_{\mathcal{V}}H_{n}$ is the classifying
map. 

The left square commutes by the definition of the \emph{fibered} product,
$G_{n}\boxtimes^{L_{\mathcal{V}}\forgetful{H_{n}}}L_{\mathcal{V}}\bar{H}_{n}$
--- see definition~\ref{def:catfiberedproduct}. The right square
commutes by the naturality of cogenerating maps. 

Now, note that the composite $H_{n}\xrightarrow{\alpha}L_{\mathcal{V}}H_{n}\xrightarrow{\varepsilon}H_{n}$
is the \emph{identity} map (a universal property of classifying maps
of coalgebras). It follows that, \emph{as a chain-map,} the composite\[
G_{n}\boxtimes^{L_{\mathcal{V}}\forgetful{H_{n}}}L_{\mathcal{V}}\bar{H}_{n}\xrightarrow{p_{n+1}}G_{n}\to H_{n}\]
 coincides with a chain-map that factors through the \emph{contractible}
chain-complex $\bar{H}_{n}$.
\end{proof}
Our main result is:
\begin{lem}
\label{lem:cofreeresolution} Let $f:C\to D$ be a cofibration as
in definition~\ref{def:rug-resolution} with rug-resolution $\{f_{i}\}:\underline{\{C\}}\to\{G_{i}\}\to D$.
Then\[
f_{\infty}=\ilimit f_{n}:C\to\ilimit G_{n}\]
 is a trivial cofibration.\end{lem}
\begin{proof}
We make extensive use of the material in appendix~\ref{sub:Limits-and-colimits}
to show that the cofiber of\[
f_{\infty}:C\to\ilimit G_{n}\]
is contractible. We focus on the category $\coalgcat$ --- the argument
in $\ircoalgcat$ is very similar. In this case, the cofiber is simply
the quotient. We will consistently use the notation $\bar{H_{n}}=\Sigma^{-1}\cone{\forgetful{H_{n}}}$

First, note that the maps\[
G_{n+1}\to G_{n}\]
 induce compatible maps\begin{eqnarray*}
L_{\mathcal{V}}\forgetful{H_{n+1}} & \to & L_{\mathcal{V}}\forgetful{H_{n}}\\
L_{\mathcal{V}}\bar{H}_{n+1} & \to & L_{\mathcal{V}}\bar{H}_{n}\end{eqnarray*}
 so proposition~\ref{pro:ilimitfibereddia} implies that\[
\ilimit G_{n}=(\ilimit G_{n})\boxtimes^{(\ilimit L_{\mathcal{V}}\forgetful{H_{n}})}(\ilimit L_{\mathcal{V}}\bar{H}_{n})\]
 and theorem~\ref{th:cofreelimits} implies that\begin{eqnarray*}
\ilimit L_{\mathcal{V}}\forgetful{H_{n}} & = & L_{\mathcal{V}}(\ilimit\forgetful{H_{i}})\\
\ilimit L_{\mathcal{V}}\bar{H}_{n} & = & L_{\mathcal{V}}(\ilimit\bar{H}_{i})=L_{\mathcal{V}}(\Sigma^{-1}\cone{\ilimit\forgetful{H_{i}}})\end{eqnarray*}
from which we conclude\[
\ilimit G_{n}=(\ilimit G_{n})\boxtimes^{L_{\mathcal{V}}(\ilimit\forgetful{H_{i}})}L_{\mathcal{V}}(\Sigma^{-1}\cone{\ilimit\forgetful{H_{i}}})\]

We claim that the projection\begin{equation}
h:(\ilimit G_{n})\boxtimes^{L_{\mathcal{V}}(\ilimit\forgetful{H_{i}})}L_{\mathcal{V}}(\Sigma^{-1}\cone{\ilimit\forgetful{H_{i}}})\to\ilimit G_{n}\label{eq:initialresolutionprojection}\end{equation}
is \emph{split} by a coalgebra morphism. To see this, first note that,
by proposition~\ref{pro:hnnummhomotopic}, each of the maps\[
H_{n+1}\twoheadrightarrow H_{n}\]
is nullhomotopic via a nullhomotopy compatible with the maps in the
inverse system $\{H_{n}\}$. This implies that \[
\ilimit\forgetful{H_{n}}\]
--- the inverse limit of \emph{chain complexes} --- is contractible.
It follows that the projection\[
\Sigma^{-1}\cone{\ilimit\forgetful{H_{n}}}\twoheadrightarrow\ilimit\forgetful{H_{n}}\]
is a trivial fibration in $\chaincat$, hence split by a map\begin{equation}
j:\ilimit\forgetful{H_{n}}\to\Sigma^{-1}\cone{\ilimit\forgetful{H_{n}}}\label{eq:splitchaincomplexes}\end{equation}
 This, in turn, induces a coalgebra morphism\[
L_{\mathcal{V}}j:L_{\mathcal{V}}(\ilimit\forgetful{H_{n}})\to L_{\mathcal{V}}(\Sigma^{-1}\cone{\ilimit\forgetful{H_{n}}})\]
 splitting the canonical surjection\[
L_{\mathcal{V}}(\Sigma^{-1}\cone{\ilimit\forgetful{H_{n}}})\twoheadrightarrow L_{\mathcal{V}}(\ilimit\forgetful{H_{n}})\]
and induces a map, $g$\begin{multline}
\ilimit G_{n}=(\ilimit G_{n})\boxtimes^{L_{\mathcal{V}}(\ilimit\forgetful{H_{i}})}L_{\mathcal{V}}(\ilimit\forgetful{H_{i}})\\
\xrightarrow{1\boxtimes L_{\mathcal{V}}j}(\ilimit G_{n})\boxtimes^{L_{\mathcal{V}}(\ilimit\forgetful{H_{i}})}L_{\mathcal{V}}(\Sigma^{-1}\cone{\ilimit\forgetful{H_{i}}})\label{eq:inducedsplitcoalgebras}\end{multline}
splitting the projection in formula~\ref{eq:initialresolutionprojection}.
Since the image of $f_{\infty}(C)$ vanishes in $L_{\mathcal{V}}(\ilimit\forgetful{H_{i}})$,
it is not hard to see that $1\boxtimes L_{\mathcal{V}}j$ is compatible
with the inclusion of $C$ in $\ilimit G_{i}$.

Now consider the diagram \begingroup\small\xymatrixcolsep{8pt}\[
\xymatrix{{(\ilimit G_{n})/f_{\infty}(C)}\ar[r]^{g\qquad\qquad\qquad\qquad\qquad}\ar@{=}[dddd] & {\left((\ilimit G_{n})\boxtimes^{L_{\mathcal{V}}(\ilimit\forgetful{H_{i}})}L_{\mathcal{V}}(\Sigma^{-1}\cone{\ilimit\forgetful{H_{i}}})\right)/f_{\infty}(C)}\ar[d]^{q}\\
{} & {(\ilimit G_{n})/f_{\infty}(C)\boxtimes^{L_{\mathcal{V}}(\ilimit\forgetful{H_{i}})}L_{\mathcal{V}}(\Sigma^{-1}\cone{\ilimit\forgetful{H_{i}}})}\ar@{=}[d]\\
{} & {(\ilimit G_{n}/f_{n}(C))\boxtimes^{L_{\mathcal{V}}(\ilimit\forgetful{H_{i}})}L_{\mathcal{V}}(\Sigma^{-1}\cone{\ilimit\forgetful{H_{i}}})}\ar@{=}[d]\\
{} & {(\ilimit H_{n})\boxtimes^{L_{\mathcal{V}}(\ilimit\forgetful{H_{i}})}L_{\mathcal{V}}(\Sigma^{-1}\cone{\ilimit\forgetful{H_{i}}})}\ar[d]^{p}\\
{(\ilimit G_{n})/f_{\infty}(C)}\ar@{=}[r] & {\ilimit H_{n}}}
\]
 \endgroup where:
\begin{enumerate}
\item The map\begin{multline*}
q:\left((\ilimit G_{n})\boxtimes^{L_{\mathcal{V}}(\ilimit\forgetful{H_{i}})}L_{\mathcal{V}}(\Sigma^{-1}\cone{\ilimit\forgetful{H_{i}}})\right)/f_{\infty}(C)\\
\to(\ilimit G_{n})/f_{\infty}(C)\boxtimes^{L_{\mathcal{V}}(\ilimit\forgetful{H_{i}})}L_{\mathcal{V}}(\Sigma^{-1}\cone{\ilimit\forgetful{H_{i}}})\end{multline*}
is induced by the projections \begingroup\xymatrixcolsep={-60pt} 
\[
\xymatrix{{} & {\ilimit G_{n}}\\
{(\ilimit G_{n})\boxtimes^{L_{\mathcal{V}}(\ilimit\forgetful{H_{i}})}L_{\mathcal{V}}(\Sigma^{-1}\cone{\ilimit\forgetful{H_{i}}})}\ar[ru]\ar[rd] & {}\\
{} & {L_{\mathcal{V}}(\Sigma^{-1}\cone{\ilimit\forgetful{H_{i}}})}}
\]
\endgroup the fact that the image of $f_{\infty}$ is effectively
only in the factor $\ilimit G_{n}$, and the defining property of
fibered products.
\item The equivalence\[
\ilimit G_{n}/f_{\infty}(C)=\ilimit G_{n}/f_{n}(C)\]
follows from theorem~\ref{thm:aequivforgetfuldirectsummand}.
\item The vertical map on the left is the identity map because $g$ splits
the map $h$ in formula~\ref{eq:initialresolutionprojection}.
\end{enumerate}
We claim that the map (projection to the left factor) \[
\forgetful p:\forgetful{(\ilimit H_{n})\boxtimes^{L_{\mathcal{V}}(\ilimit\forgetful{H_{i}})}L_{\mathcal{V}}(\Sigma^{-1}\cone{\ilimit\forgetful{H_{i}}})}\to\forgetful{\ilimit H_{n}}\]
is \emph{nullhomotopic} (as a $\chaincat$-morphism). This follows
immediately from the fact that \[
\ilimit H_{n}\hookrightarrow L_{\mathcal{V}}(\ilimit\forgetful{H_{n}})\]
by corollary~\ref{cor:inverselimitintersectiontwo}, so that\begin{multline*}
(\ilimit H_{n})\boxtimes^{L_{\mathcal{V}}(\ilimit\forgetful{H_{i}})}L_{\mathcal{V}}(\Sigma^{-1}\cone{\ilimit\forgetful{H_{i}}})\\
\subseteq L_{\mathcal{V}}(\ilimit\forgetful{H_{i}})\boxtimes^{L_{\mathcal{V}}(\ilimit\forgetful{H_{i}})}L_{\mathcal{V}}(\Sigma^{-1}\cone{\ilimit\forgetful{H_{i}}})\\
=L_{\mathcal{V}}(\Sigma^{-1}\cone{\ilimit\forgetful{H_{i}}})\end{multline*}
and $L_{\mathcal{V}}(\Sigma^{-1}\cone{\ilimit\forgetful{H_{i}}})$
is contractible, by proposition~\ref{prop:lefthomotopy} and the
contractibility of $\Sigma^{-1}\cone{\ilimit\forgetful{H_{i}}}$.

We conclude that \[
(\ilimit G_{n})/f_{\infty}(C)\xrightarrow{\mathrm{id}}(\ilimit G_{n})/f_{\infty}(C)\]
is nullhomotopic so $(\ilimit G_{n})/f_{\infty}(C)$ is contractible
and\[
\forgetful{f_{\infty}}:\forgetful C\to\forgetful{\ilimit G_{n}}\]
 is a weak equivalence in $\chaincat$, hence (by definition~\ref{def:cofibrationweakequiv})
$f_{\infty}$ is a weak equivalence.\end{proof}
\begin{cor}
\label{cor:factortrivialcofibrationfibration}Let $\mathcal{V}=\{\mathcal{V}(n)\}$
be a $\Sigma$-cofibrant operad such that $\mathcal{V}(n)$ is of
finite type for all $n\ge0$. Let\[
f:A\to B\]
be a morphism in $\ircoalgcat$ or $\coalgcat$. Then there exists
a functorial factorization of $f$ \[
A\to Z(f)\to B\]
where \[
A\to Z(f)\]
 is a trivial cofibration and \[
Z(f)\to B\]
is a fibration. \end{cor}
\begin{rem*}
This is condition CM5, statement 1 in the definition of a model category
at the beginning of this section. It, therefore, proves that the model
structure described in \ref{def:cofibrationweakequiv} and \ref{def:fibration}
is well-defined.

By abuse of notation, we will call the $\{f_{i}\}:\underline{\{A\}}\to\{G_{i}\}\to L_{\mathcal{V}}(\forgetful A)\boxtimes B\to B$
the \emph{rug-resolution of the morphism} $A\to B$ (see the proof
below), where $\{f_{i}\}:\underline{\{A\}}\to\{G_{i}\}\to L_{\mathcal{V}}(\forgetful A)\boxtimes B$
is the rug-resolution of the cofibraton $A\to L_{\mathcal{V}}(\forgetful A)\boxtimes B$.

See proposition~\ref{th:cofreelimits} and corollary~\ref{cor:generalcoalglimits}
for the definition of inverse limit in the category $\mathrf{C}$.\end{rem*}
\begin{proof}
Simply apply definition~\ref{def:rug-resolution}and lemma~\ref{lem:cofreeresolution}
to the cofibration\[
A\to L_{\mathcal{V}}(\forgetful A)\boxtimes B\]
and project to the second factor.
\end{proof}
We can characterize fibrations now:
\begin{cor}
\label{cor:fibrationcharacterization}If $\mathcal{V}=\{\mathcal{V}(n)\}$
is a $\Sigma$-cofibrant operad such that $\mathcal{V}(n)$ is of
finite type for all $n\ge0$, then all fibrations are retracts of
their rug-resolutions.\end{cor}
\begin{rem*}
This shows that rug-resolutions of maps contain canonical fibrations
and all others are retracts of them.\end{rem*}
\begin{proof}
Suppose $p:A\to B$ is some fibration. We apply corollary~\ref{cor:factortrivialcofibrationfibration}
to it to get a commutative diagram \[
\xymatrix{{A}\ar[d]_{i} & {}\\
{\bar{A}}\ar[r]_{a_{\infty}} & {B}}
\]
 where $i:A\to\bar{A}$ is a trivial cofibration and $u_{\infty}:\bar{A}\to B$
is a fibration. We can complete this to get the diagram \[
\xymatrix{{A}\ar[d]_{i}\ar@{=}[r] & {A}\ar[d]^{p}\\
{\bar{A}}\ar[r]_{a_{\infty}}\ar@{.>}[ru] & {B}}
\]
The fact that $p:A\to B$ is a fibration and definition~\ref{def:fibration}
imply the existence of the dotted arrow making the whole diagram commute.
But this \emph{splits} the inclusion $i:A\to\bar{A}$ and implies
the result.
\end{proof}
The rest of this section will be spent on \emph{trivial} fibrations
--- with a mind to proving the \emph{second statement} in CM 4 in
theorem~\ref{thm:secondcm4}. Recall that the first statement was
a consequence of our \emph{definition} of fibrations in $\coalgcat$
and $\ircoalgcat$.

\subsection{Proof of CM~4\label{sub:Proof-of-CM4}}

The first part of CM~4 is trivial: we have \emph{defined} fibrations
as morphisms that satisfy it --- see definition~\ref{def:fibration}.
The proof of the second statement of CM~4 is more difficult and makes
extensive use of the Rug Resolution defined in definition~\ref{def:rug-resolution}.

We begin by showing that a fibration of coalgebras becomes a fibration
in $\chaincat$ under the forgetful functor:
\begin{prop}
\label{pro:coalgfibrationchaincatfibration}Let $p:A\to B$ be a fibration
in $\mathrf{C}=\ircoalgcat$ or $\coalgcat$. Then \[
\forgetful p:\forgetful A\to\forgetful B\]
is a fibration in $\chaincat^{+}$ or $\chaincat$, respectively.\end{prop}
\begin{proof}
In the light of corollary~\ref{cor:fibrationcharacterization}, it
suffices to prove this for rug-resolutions of fibrations. 

Since they are iterated pullbacks of fibrations with contractible
total spaces, it suffices to prove the result for something of the
form\[
A\boxtimes^{L_{\mathcal{V}}B}L_{\mathcal{V}}(\Sigma^{-1}\cone B)\to A\]
where $f:A\to L_{\mathcal{V}}B$ is some morphism. The fact that all
morphisms are coalgebra morphisms implies the existence of a coalgebra
structure on \[
Z=\forgetful A\oplus^{\forgetful{L_{\mathcal{V}}B}}\forgetful{L_{\mathcal{V}}(\Sigma^{-1}\cone B)}\subset\forgetful{A\boxtimes L_{\mathcal{V}}(\Sigma^{-1}\cone B)}\rightrightarrows\forgetful{L_{\mathcal{V}}B}\]
where $\forgetful A\oplus^{\forgetful{L_{\mathcal{V}}B}}\forgetful{L_{\mathcal{V}}(\Sigma^{-1}\cone B)}$
is the fibered product in $\chaincat$. Since $L_{\mathcal{V}}(\Sigma^{-1}\cone B)\to L_{\mathcal{V}}B$
is surjective, (because it is induced by the surjection, $\Sigma\cone B)\to B$)
it follows that the equalizer \[
\forgetful{A\boxtimes L_{\mathcal{V}}(\Sigma^{-1}\cone B)}\rightrightarrows\forgetful{L_{\mathcal{V}}B}\]
surjects onto $\forgetful A$. Since $Z$ has a coalgebra structure,
it is contained in the core, \[
\core{\forgetful{A\boxtimes L_{\mathcal{V}}(\Sigma^{-1}\cone B)}\rightrightarrows\forgetful{L_{\mathcal{V}}B}}=A\boxtimes^{L_{\mathcal{V}}B}L_{\mathcal{V}}(\Sigma^{-1}\cone B)\]
which also surjects onto $A$ --- so the projection \[
\forgetful{A\boxtimes^{L_{\mathcal{V}}B}L_{\mathcal{V}}(\Sigma^{-1}\cone B)}\to\forgetful A\]
 is surjective and --- as a map of graded $\ring$-modules --- split.
This is the definition of a fibration in $\chaincat$.
\end{proof}
We are now in a position to prove the second part of CM~4:
\begin{thm}
\label{thm:secondcm4}Given a commutative solid arrow diagram \[
\xymatrix{{U}\ar[r]^{f}\ar[d]_{i} & {A}\ar[d]^{p}\\
{W}\ar[r]_{g}\ar@{.>}[ur] & {B}}
\]
 where $i$ is a any cofibration and $p$ is a \emph{trivial} fibration,
the dotted arrow exists.\end{thm}
\begin{proof}
Because of corollary~\ref{cor:fibrationcharacterization}, it suffices
to prove the result for the \emph{rug-resolution} of the trivial fibration
$p:A\to B$. We begin by considering the diagram\[
\xymatrix{{\forgetful U}\ar[r]^{\forgetful f}\ar[d]_{\forgetful i} & {\forgetful A}\ar[d]^{\forgetful p}\\
{\forgetful W}\ar[r]_{\forgetful g}\ar@{.>}[ur]^{\ell} & {\forgetful B}}
\]
Because of proposition~\ref{pro:coalgfibrationchaincatfibration},
$\forgetful p$ is a trivial fibration and the dotted arrow exists
in $\chaincat$. 

If $\alpha:A\to L_{\mathcal{V}}A$ is the classifying map of $A$,
$\hat{\ell}:W\to L_{\mathcal{V}}A$ is induced by $\ell:\forgetful W\to\forgetful A$,
and $p_{2}:L_{\mathcal{V}}A\boxtimes B\to B$ is projection to the
second factor, we get a commutative diagram\begin{equation}
\xymatrix{\xyC{70pt}{U}\ar[r]^{(\alpha\circ f)\boxtimes(p\circ f)}\ar[d]_{i} & {L_{\mathcal{V}}A\boxtimes B}\ar[d]^{p_{2}}\\
{W}\ar[r]_{g}\ar@{.>}[ur]_{\hat{\ell}\boxtimes g} & {B}}
\label{eq:cofibrationdia1}\end{equation}

It will be useful to build the rug-resolutions of $A\to L_{\mathcal{V}}A\boxtimes B=G_{0}$
and $B\to L_{\mathcal{V}}B\boxtimes B=\tilde{G}_{0}$ in parallel
--- denoted $\{G_{n}\}$ and $\{\tilde{G}_{n}\}$, respectively. Clearly
the vertical morphisms in \begin{equation}
\xymatrix{{A}\ar[r]^{f_{0}\qquad}\ar[d]_{p} & {L_{\mathcal{V}}A\boxtimes B}\ar[d]^{q_{0}}\\
{B}\ar[r]_{\bar{f}_{0}\qquad} & {L_{\mathcal{V}}B\boxtimes B}}
\label{eq:startinduction}\end{equation}
are trivial fibrations via strict homotopies --- see propositions~\ref{prop:lefthomotopy}
and \ref{pro:stricthomotopyproducts}.

We prove the result by an induction that:
\begin{enumerate}
\item lifts the map $\hat{\ell}_{1}=\hat{\ell}\boxtimes g:W\to L_{\mathcal{V}}A\boxtimes B$
to norphisms $\hat{\ell}_{k}:W\to G_{k}$ to successively higher stages
of the rug-resolution of $p$. Diagram~\ref{eq:startinduction} implies
the base case.
\item establishes that the vertical morphisms in\[
\xymatrix{{A}\ar[r]^{f_{n}}\ar[d]_{p} & {G_{n}}\ar[d]^{q_{n}}\\
{B}\ar[r]_{\bar{f}_{n}} & {\tilde{G}_{n}}}
\]
are trivial fibrations for all $n$, via \emph{strict} homotopies. 
\end{enumerate}
Lemma~\ref{lem:reltrivfibrationsplit}implies that we can find splitting
maps $u$ and $v$ such that \begin{equation}
\xymatrix{{\forgetful A}\ar[r]^{f_{n}} & {\forgetful{G_{n}}}\\
{\forgetful B}\ar[r]_{\bar{f}_{n}}\ar[u]^{u} & {\forgetful{\tilde{G}_{n}}}\ar[u]_{v}}
\label{eq:goodsplittings}\end{equation}
$p\circ u=1:\forgetful B\to\forgetful B$, $q_{n}\circ v:\forgetful{\bar{G}_{n}}\to\forgetful{\bar{G}_{n}}$
and contracting homotopies $\Phi_{1}$ and $\Phi_{2}$ such that\begin{equation}
\xymatrix{{\forgetful A\otimes I}\ar[r]^{f_{n}\otimes1}\ar[d]_{\Phi_{1}} & {\forgetful{G_{n}}\otimes I}\ar[d]^{\Phi_{2}}\\
{\forgetful A}\ar[r]_{f_{n}} & {\forgetful{G_{n}}}}
\label{eq:goodhomotopies}\end{equation}
commutes, where $d\Phi_{1}=u\circ p-1$, and $d\Phi_{2}=v\circ q_{n}-1$
--- where $\Phi_{1}$ can be specified beforehand. Forming quotients
gives rise to a commutative diagram\[
\xymatrix{{\forgetful A}\ar[r]^{f_{n}}\ar[d]_{p} & {\forgetful{G_{n}}}\ar[d]^{q_{n}}\ar[r] & H_{n}\ar[d]^{\hat{q}_{n}}\\
{\forgetful B}\ar[r]_{\bar{f}_{n}} & {\forgetful{\tilde{G}_{n}}}\ar[r] & \tilde{H}_{n}}
\]
Furthermore the commutativity of diagrams~\ref{eq:goodsplittings}
and \ref{eq:goodhomotopies} implies that $v$ induces a splitting
map $w:\tilde{H}_{n}\to H_{n}$ and $\Phi_{2}$ induces a homotopy
$\Xi:H_{n}\otimes I\to H_{n}$ with $d\Xi=v\circ\hat{q}_{n}-1$ ---
so $\hat{q}_{n}$ is a \emph{weak equivalence} in $\chaincat$ ---
even a trivial fibration.

If we assume that the lifting has been carried out to the $n^{\text{th}}$
stage, we have a map\[
\ell_{n}:W\to G_{n}\]
making \[
\xymatrix{{\forgetful W}\ar[r]^{\ell_{n}}\ar[d]_{g} & {\forgetful{G_{n}}}\ar[d]^{q_{n}}\ar[r] & H_{n}\ar[d]^{\hat{q}_{n}}\\
{\forgetful B}\ar[r]_{\bar{f}_{n}} & {\forgetful{\tilde{G}_{n}}}\ar[r] & \tilde{H}_{n}}
\]
commute. Since the image of $B$ in $\tilde{H}_{n}$ vanishes (by
the way $H_{n}$ and $\tilde{H}_{n}$ are constructed --- see statement
2 in definition~\ref{def:rug-resolution}), it follows that the image
of $W$ in $H_{n}$ lies in the \emph{kernel} of $\hat{q}_{n}$ ---
a trivial fibration in $\chaincat$. We conclude that the inclusion
of $W$ in $H_{n}$ is null-homotopic, hence \emph{lifts} to $\bar{H}_{n}=\Sigma^{-1}\cone{H_{n}}$
in such a way that \[
\xymatrix{ &  & \bar{H}_{n}\ar[d]\\
{\forgetful W}\ar[r]^{\ell_{n}}\ar[d]_{g}\ar@/^{9pt}/[rru]^{r} & {\forgetful{G_{n}}}\ar[d]^{q_{n}}\ar[r]^{t} & H_{n}\ar[d]^{\hat{q}_{n}}\\
{\forgetful B}\ar[r]_{\bar{f}_{n}} & {\forgetful{\tilde{G}_{n}}}\ar[r] & \tilde{H}_{n}}
\]
commutes --- as a diagram of \emph{chain-complexes}. Now note that
$G_{n+1}$ is the fibered product $G_{n}\boxtimes^{L_{\mathcal{V}}(G_{n}/A)}L_{\mathcal{V}}\Sigma^{-1}\cone{G_{n}/A}$
--- and that the chain-maps $r$ and $t\circ\ell_{n}$ induce c\emph{oalgebra
morphisms} making the diagram\[
\xymatrix{ & \bar{H}_{n}\ar[d]\\
{W}\ar[rd]_{\ell_{n}}\ar[ru]^{L_{\mathcal{V}}r} & H_{n}\\
 & {G_{n}}\ar[u]_{L_{\mathcal{V}}t}}
\]
commute --- thereby inducing a coagebra-morphism\[
\ell_{n+1}:W\to G_{n}\boxtimes^{L_{\mathcal{V}}(H_{n})}L_{\mathcal{V}}\Sigma^{-1}\cone{H_{n}}=G_{n+1}\]
that makes the diagram\[
\xymatrix{ & G_{n+1}\ar[d]^{p_{n+1}}\\
{W}\ar[r]_{\ell_{n}}\ar[ru]^{\ell_{n+1}} & G_{n}}
\]
commute (see statement 5 of definition~\ref{def:rug-resolution}).
This proves assertion 1 in the induction step. 

To prove assertion 2 in our induction hypothesis, note that the natural
homotopy in diagram~\ref{eq:goodhomotopies} induces (by passage
to the quotient) a natural homotopy, $\Phi'$, that makes the diagram
of chain-complexes\[
\xymatrix{{\forgetful{G_{n}}\otimes I}\ar[d]_{\Phi_{2}}\ar[r]^{t\otimes1} & H_{n}\otimes I\ar[d]^{\Phi'}\\
{\forgetful{G_{n}}}\ar[r]_{t} & H_{n}}
\]
commute. This can be expanded to a commutative diagram \[
\xymatrix{{\forgetful{G_{n}}\otimes I}\ar[d]_{\Phi_{2}}\ar[r]^{t\otimes1} & H_{n}\otimes I\ar[d]^{\Phi'} & \bar{H}_{n}\otimes I\ar[d]^{\bar{\Phi}}\ar[l]\\
{\forgetful{G_{n}}}\ar[r]_{t} & H_{n} & \bar{H}_{n}\ar[l]}
\]
The conclusion follows from the fact that $\Phi'$ and $\bar{\Phi}$
induce \emph{strict} homotopies (see definition~\ref{def:strictequivalence})
after the cofree coalgebra-functor is applied (see proposition~\ref{prop:lefthomotopy})
and proposition~\ref{pro:strictfiberedprodinvar}.

We conclude that $G_{n+1}\to\tilde{G}_{n+1}$ is a trivial fibration.

Induction shows that we can define a lifting\[
\ell_{\infty}:W\to\ilimit G_{n}\]
that makes the diagram\[
\xymatrix{\xyC{70pt}{U}\ar[r]^{\iota\circ f}\ar[d]_{i} & {\ilimit G_{n}}\ar[d]^{p_{\infty}}\\
{W}\ar[r]_{g}\ar@{.>}[ur]_{\ell_{\infty}} & {B}}
\]
commute.
\end{proof}

\subsection{The bounded case\label{sec:The-bounded-case}}

In this section, we develop a model structure on a category of coalgebras
whose underlying chain-complexes are b\emph{ounded from below.} 
\begin{defn}
\label{def:boundedcategories}Let:
\begin{enumerate}
\item $\bchaincat$ denote the subcategory of $\chaincat$ bounded at dimension
$1$. If $A\in\bchaincat$, then $A_{i}=0$ for $i<1$.
\item $\bircoalgcat$ denote the category of pointed irreducible coalgebras,
$C$, over $\mathcal{V}$ such that $\forgetful C\in\bchaincat$.
This means that $C_{i}=0$, $i<1$. Note, from the definition of $\forgetful C$
as the kernel of the augmentation map, that the underlying chain-complex
of $C$ is equal to $\ring$ in dimension $0$.
\end{enumerate}
\end{defn}
There is clearly an inclusion of categories\[
\iota:\bchaincat\to\chaincat\]
compatible with model structures.

Now we define our model structure on $\bircoalgcat$:
\begin{defn}
\label{def:bcofibrationweakequiv}A morphism $f:A\to B$ in $\bircoalgcat$
will be called
\begin{enumerate}
\item a \emph{weak equivalence} if $\forgetful f:\forgetful A\to\forgetful B$
is a weak equivalence in $\bchaincat$ (i.e., a chain homotopy equivalence).
An object $A$ will be called \emph{contractible} if the augmentation
map\[
A\to\ring\]
is a weak equivalence. 
\item a \emph{cofibration} if $\forgetful f$ is a cofibration in $\bchaincat$.
\item a \emph{trivial cofibration} if it is a weak equivalence and a cofibration.
\end{enumerate}
\end{defn}
\begin{rem*}
A morphism is a cofibration if it is a degreewise split monomorphism
of chain-complexes. Note that all objects of $\bircoalgcat$ are cofibrant.

If $\ring$ is a \emph{field}, all modules are vector spaces therefore
free. Homology equivalences of bounded free chain-complexes induce
chain-homotopy equivalence, so our notion of weak equivalence becomes
the \emph{same} as homology equivalence (or quasi-isomorphism).\end{rem*}
\begin{defn}
\label{def:bfibration}A morphism $f:A\to B$ in $\bircoalgcat$ will
be called 
\begin{enumerate}
\item a \emph{fibration} if the dotted arrow exists in every diagram of
the form \[
\xymatrix{{U}\ar[r]\ar[d]_{i} & {A}\ar[d]^{f}\\
{W}\ar[r]\ar@{.>}[ur] & {B}}
\]
 in which $i:U\to W$ is a trivial cofibration.
\item a \emph{trivial fibration} if it is a fibration and a weak equivalence.
\end{enumerate}
\end{defn}
\begin{cor}
\label{cor:boundedproof}If $\mathcal{V}=\{\mathcal{V}(n)\}$ is an
operad satisfying condition~\ref{cond:mainoperadassumption}, the
description of cofibrations, fibrations, and weak equivalences given
in definitions~\ref{def:bcofibrationweakequiv} and \ref{def:bfibration}
satisfy the axioms for a model structure on $\bircoalgcat$.\end{cor}
\begin{proof}
We carry out all of the constructions of \S~\ref{sec:The-general-case}
and appendix~\ref{sec:Basic-constructions} while consistently replacing
cofree coalgebras by their truncated versions (see \cite{Smith-cofree}).
This involves substituting $M_{\mathcal{V}}(*)$ for $L_{\mathcal{V}}(*)$
and $\mathrf{F}_{\mathcal{V}}(*)$ for $P_{\mathcal{V}}(*)$ . 
\end{proof}

\section{Examples\label{sec:Examples}}

We will give a few examples of the model structure developed here.
In all cases, we will make the simplifying assumption that $\ring$
is a field (this is \emph{not} to say that interesting applications
\emph{only} occur when $\ring$ is a field). We begin with coassociative
coalgebras over the rationals:
\begin{example}
\label{exa:ircoalgoverq}Let $\mathcal{V}$ be the operad with component
$n$ equal to $\mathbb{Q}S_{n}$ with the obvious $S_{n}$-action
--- and we consider the category of pointed, irreducible coagebras,
$\ircoalgcat$. Coalgebras over this $\mathcal{V}$ are coassociative
coalgebras. In this case $P_{\mathcal{V}}C=T(C)$, the graded tensor
algebra with coproduct\begin{equation}
c_{1}\otimes\cdots\otimes c_{n}\mapsto\sum_{k=0}^{n}(c_{1}\otimes\cdots\otimes c_{k})\otimes(c_{k+1}\otimes\cdots\otimes c_{n})\label{eq:tensorcoproduct}\end{equation}
where $c_{1}\otimes\cdots\otimes c_{0}=c_{n+1}\otimes\cdots\otimes c_{n}=1\in C^{0}=\mathbb{Q}$.
The $n$-fold coproducts are just composites of this $2$-fold coproduct
and the {}``higher'' coproducts vanish identically. We claim that
this makes \begin{equation}
A\boxtimes B=A\otimes B\label{eq:aboxtimesbtensor}\end{equation}
 This is due to the well-known identity $T(\forgetful A\oplus\forgetful B)=T(\forgetful A)\otimes T(\forgetful B)$.
The category $\bircoalgcat$ is a category of 1-connected coassociative
coalgebras where weak equivalence is equivalent to homology equivalence.
\end{example}
If we assume coalgebras to be cocommutative we get:
\begin{example}
\label{exa:quillen}Suppose $\ring=\mathbb{Q}$ and $\mathcal{V}$
is the operad with all components equal to $\mathbb{Q}$, concentrated
in dimension $0$, and equipped with trivial symmetric group actions.
Coalgebras over $\mathcal{V}$ are just cocommutative, coassociative
coalgebras and $\bircoalgcat$ is a category of 1-connected coalgebras
similar to the one Quillen considered in \cite{Quillen:1969}. Consequently,
our model structure for $\bircoalgcat$ induces the model structure
defined by Quillen in \cite{Quillen:1969} on the subcategory of 2-connected
coalgebras.

In this case, $P_{\mathcal{V}}C$ is defined by\[
P_{\mathcal{V}}C=\bigoplus_{n\ge0}(C^{\otimes n})^{S_{n}}\]
where $(C^{\otimes n})^{S_{n}}$ is the submodule of \[
\underbrace{C\otimes\cdots\otimes C}_{n\text{\,\ factors}}\]
invariant under the $S_{n}$-action. The assumption that the base-ring
is $\mathbb{Q}$ implies a canonical isomorphism\[
P_{\mathcal{V}}C=\bigoplus_{n\ge0}(C^{\otimes n})^{S_{n}}\cong S(C)\]

Since $S(\forgetful A\oplus\forgetful B)\cong S(\forgetful A)\otimes S(\forgetful B)$,
we again get $A\boxtimes B=A\otimes B$. 
\end{example}
\appendix

\section{Nearly free modules\label{sec:nearlyfree}}

In this section, we will explore the class of nearly free $\integers$-modules
--- see definition~\ref{def:nearlyfree}. We show that this is closed
under the operations of taking direct sums, tensor products, countable
products and cofree coalgebras. It appears to be fairly large, then,
and it would be interesting to have a direct algebraic characterization.
\begin{rem}
\label{rem:nearlyfreeimpliesflat}A module must be torsion-free (hence
flat) to be nearly free. The converse is not true, however: $\mathbb{Q}$
is flat but \emph{not} nearly free.
\end{rem}
The definition immediately implies that:

\begingroup\def\ring{\integers}
\begin{prop}
\label{prop:nfsubmodule}Any submodule of a nearly free module is
nearly free.
\end{prop}
Nearly free modules are closed under operations that preserve free
modules: 
\begin{prop}
\label{prop:nfsumotimes}Let $M$ and $N$ be $\integers$-modules.
If they are nearly free, then so are $M\oplus N$ and $M\otimes N$.

Infinite direct sums of nearly free modules are nearly free.\end{prop}
\begin{proof}
If $F\subseteq M\oplus N$ is countable, so are its projections to
$M$ and $N$, which are free by hypothesis. It follows that $F$
is a countable submodule of a free module.

The case where $F\subseteq M\otimes N$ follows by a similar argument:
The elements of $F$ are finite linear combinations of monomials $\{m_{\alpha}\otimes n_{\alpha}\}$
--- the set of which is countable. Let\begin{eqnarray*}
A & \subseteq & M\\
B & \subseteq & N\end{eqnarray*}
 be the submodules generated, respectively, by the $\{m_{\alpha}\}$
and $\{n_{\alpha}\}$. These will be countable modules, hence $\integers$-free.
It follows that \[
F\subseteq A\otimes B\]
is a free module.

Similar reasoning proves the last statement, using the fact that any
direct sum of free modules is free.\end{proof}
\begin{prop}
\label{prop:nfprodfree}Let $\{F_{n}\}$ be a countable collection
of $\integers$-free modules. Then\[
\prod_{n=1}^{\infty}F_{n}\]
is nearly free.\end{prop}
\begin{proof}
In the case where $F_{n}=\integers$ for all $n$\[
B=\prod_{n=1}^{\infty}\integers\]
is the Baer-Specker group, which is well-known to be nearly free ---
see \cite{Baer:1937}, \cite[vol. 1, p. 94 Theorem 19.2]{Fuchs:1970},
and\cite{Blass-Gobel:1996}. It is also well-known \emph{not} to be
$\integers$-free --- see \cite{Baer-Specker-nonfree} or the survey
\cite{Baer-Specker-survey}.

First suppose each of the $F_{n}$ are countably generated. Then\[
F_{n}\subseteq B\]
and \[
\prod F_{n}\subseteq\prod B=B\]
which is nearly-free. 

In the general case, any countable submodule, $C$, of $\prod F_{n}$
projects to a countably-generated submodule, $A_{n}$, of $F_{n}$
under all of the projections\[
\prod F_{n}\to F_{n}\]
and, so is contained in \[
\prod A_{n}\]
which is nearly free, so $C$ must be $\integers$-free.\end{proof}
\begin{cor}
\label{cor:nfprodnf}Let $\{N_{k}\}$ be a countable set of nearly
free modules. Then\[
\prod_{k=1}^{\infty}N_{k}\]
is also nearly free.\end{cor}
\begin{proof}
Let\[
F\subset\prod_{k=1}^{\infty}N_{k}\]
 be countable. If $F_{k}$ is its projection to factor $N_{k}$, then
$F_{k}$will be countable, hence free. It follows that\[
F\subset\prod_{k=1}^{\infty}F_{k}\]
 and the conclusion follows from proposition~\ref{prop:nfprodfree}. \end{proof}
\begin{cor}
\label{cor:nfhom1}Let $A$ be nearly free and let $F$ be $\integers$-free
of countable rank. Then\[
\homz(F,A)\]
is nearly free.\end{cor}
\begin{proof}
This follows from corollary~\ref{cor:nfprodnf} and the fact that
\[
\homz(F,A)\cong\prod_{k=1}^{\mathrm{rank}(F)}A\]
\end{proof}
\begin{cor}
\label{cor:nfprodfn}Let $\{F_{n}\}$ be a sequence of countably-generated
$\zs n$-projective modules and and let $A$ be nearly free. Then\[
\prod_{n=1}^{\infty}\homzs n(F_{n},A^{\otimes n})\]
 is nearly free.\end{cor}
\begin{proof}
This is a direct application of the results of this section and the
fact that \[
\homzs n(F_{n},A^{\otimes n})\subseteq\homz(F_{n},A^{\otimes n})\subseteq\homz(\hat{F}_{n},A^{\otimes n})\]
where $\hat{F}_{n}$ is a $\zs n$-free module of which $F_{n}$ is
a direct summand.\end{proof}
\begin{thm}
\label{thm:nfcofreenf}Let $C$ be a nearly free $\integers$-module
and let $\mathcal{V}=\{\mathcal{V}(n)\}$ be a $\Sigma$-finite operad
with $\mathcal{V}(n)$ of finite type for all $n\ge0$. Then\begin{eqnarray*}
 & \forgetful{L_{\mathcal{V}}C}\\
 & \forgetful{M_{\mathcal{V}}C}\\
 & \forgetful{P_{\mathcal{V}}C}\\
 & \forgetful{\mathrf{F}_{\mathcal{V}}C}\end{eqnarray*}
are all nearly free.\end{thm}
\begin{proof}
This follows from theorem~\ref{cor:cofreedirectlimits} which states
that all of these are submodules of \[
\prod_{n\ge0}(\mathcal{V}(n),A^{\otimes n})\]
and the fact that near-freeness is inherited by submodules.
\end{proof}
\endgroup

\section{Category-theoretic constructions\label{sec:Basic-constructions}}

In this section, we will study general properties of coalgebras over
an operad. Some of the results will require coalgebras to be pointed
irreducible. We begin by recalling the structure of cofree coalgebras
over operads in the pointed irreducible case.

\subsection{Cofree-coalgebras}

We will make extensive use of \emph{cofree coalgebras} over an operad
in this section --- see definition~\ref{def:cofreecoalgebra}.

If they exist, it is not hard to see that cofree coalgebras must be
\emph{unique} up to an isomorphism.

The paper \cite{Smith-cofree} gave an explicit construction of $L_{\mathcal{V}}C$
when $C$ was an $\ring$-free chain complex. When $\ring$ is a field,
all chain-complexes are $\ring$-free, so the results of the present
paper are already true in that case. 

Consequently, we will restrict ourselves to the case where $\ring=\integers$. 
\begin{prop}
The forgetful functor (defined in definition~\ref{def:forgetful})
and cofree coalgebra functors define adjoint pairs\begin{eqnarray*}
\pcoalg V{\ast}:\chaincat & \leftrightarrows & \ircoalgcat:\forgetful{\ast}\\
L_{\mathcal{V}}(\ast):\chaincat & \leftrightarrows & \coalgcat:\forgetful{\ast}\end{eqnarray*}
\end{prop}
\begin{rem*}
The adjointness of the functors follows from the universal property
of cofree coalgebras --- see~\cite{Smith-cofree}. 
\end{rem*}
The Adjoints and Limits Theorem in \cite{MacLane:cw} implies that:
\begin{thm}
\label{th:cofreelimits}If $\{A_{i}\}\in\mathrm{ind}-\chaincat$ and
$\{C_{i}\}\in\dircoalgcat$ or $\dcoalgcat$ then\begin{eqnarray*}
\ilimit P_{\mathcal{V}}(A_{i}) & = & P_{\mathcal{V}}(\ilimit A_{i})\\
\ilimit L_{\mathcal{V}}(A_{i}) & = & L_{\mathcal{V}}(\ilimit A_{i})\\
\forgetful{\dlimit C_{i}} & = & \dlimit\forgetful{C_{i}}\end{eqnarray*}
\end{thm}
\begin{rem*}
This implies that colimits in $\ircoalgcat$ or $\coalgcat$ are the
same as colimits of underlying chain-complexes. \end{rem*}
\begin{prop}
\label{pro:directlimitfintelygenerated}If $C\in\chaincat$, let $\mathrf{G}(C)$
denote the lattice of countable subcomplexes of $C$. Then\[
C=\dlimit\mathrf{G}(C)\]
\end{prop}
\begin{proof}
Clearly $\dlimit\mathrf{G}(C)\subseteq C$ since all of the canonical
maps to $C$ are inclusions. Equality follows from every element $x\in C_{k}$
being contained in a finitely generated subcomplex, $C_{x}$, defined
by

\[
(C_{x})_{i}=\begin{cases}
\ring\cdot x & \text{ if }i=k\\
\ring\cdot\partial x & \text{ if }i=k-1\\
0 & \text{otherwise}\end{cases}\]
\end{proof}
\begin{lem}
\label{lem:homdirectlimit}Let $n>1$ be an integer, let $F$ be a
finitely-generated projective (non-graded) $\zs n$-module, and let
$\{C_{\alpha}\}$ a direct system of modules. Then the natural map\[
\dlimit\homzs n(F,C_{\alpha})\to\homzs n(F,\dlimit C_{\alpha})\]
is an isomorphism.

If $F$ and the $\{C_{\alpha}\}$ are graded, the corresponding statement
is true if $F$ is finitely-generated and $\zs n$-projective in each
dimension.\end{lem}
\begin{proof}
We will only prove the non-graded case. The graded case follows from
the fact that the maps of the $\{C_{\alpha}\}$ preserve grade.

In the non-graded case, finite generation of $F$ implies that the
natural map\[
\bigoplus_{\alpha}\homzs n(F,C_{\alpha})\to\homzs n(F,\bigoplus_{\alpha}C_{\alpha})\]
is an isomorphism, where $\alpha$ runs over any indexing set. The
projectivity of $F$ implies that $\homzs n(F,*)$ is exact, so the
short exact sequence defining the filtered colimit is preserved.\end{proof}
\begin{prop}
\label{prop:proddirectlimit}Let $\mathcal{V}=\{\mathcal{V}(n)\}$
be an operad satisfying condition~\ref{cond:mainoperadassumption},
and let $C$ be a chain-complex with $\mathrf{G}(C)=\{C_{\alpha}\}$
a family of flat subcomplexes ordered by inclusion that is closed
under countable sums. In addition, suppose\[
C=\dlimit C_{\alpha}\]
Then\[
\prod_{n\ge0}\homzs n(\mathcal{V}(n),C^{\otimes n})=\dlimit\prod_{n\ge0}\homzs n(\mathcal{V}(n),C_{\alpha}^{\otimes n})\]
\end{prop}
\begin{proof}
Note that $C$, as the limit of flat modules, it itself flat.

The $\integers$-flatness of $C$ implies that any $y\in C^{\otimes n}$
is in the image of\[
C_{\alpha}^{\otimes n}\hookrightarrow C^{\otimes n}\]
for some $C_{\alpha}\in\mathrf{G}(C)$ and any $n\ge0$. The finite
generation and projectivity of the $\{\mathcal{V}(n)\}$ in every
dimension implies that any map\[
x_{i}\in\homzs n(\mathcal{V}(n),C^{\otimes n})_{j}\]
lies in the image of \[
\homzs n(\mathcal{V}(n),C_{\alpha_{i}}^{\otimes n})\hookrightarrow\homzs n(\mathcal{V}(n),C^{\otimes n})\]
for some $C_{\alpha_{i}}\in\mathrf{G}(C)$. This implies that\[
x\in\homzs n(\mathcal{V}(n),C^{\otimes n})\]
 lies in the image of\[
\homzs n(\mathcal{V}(n),C_{\alpha}^{\otimes n})\hookrightarrow\homzs n(\mathcal{V}(n),C^{\otimes n})\]
where $C_{\alpha}=\sum_{i=0}^{\infty}C_{\alpha_{i}}$, which is still
a subcomplex of the lattice $\mathrf{G}(C)$.

If \[
x=\prod x_{n}\in\prod_{n\ge0}\homzs n(\mathcal{V}(n),C^{\otimes n})\]
then each $x_{n}$ lies in the image of\[
\homzs n(\mathcal{V}(n),C_{\alpha_{n}}^{\otimes n})\hookrightarrow\homzs n(\mathcal{V}(n),C^{\otimes n})\]
 where $C_{\alpha_{n}}\in\mathrf{G}(C)$ and $x$ lies in the image
of \[
\prod_{n\ge0}\homzs n(\mathcal{V}(n),C_{\alpha}^{\otimes n})\hookrightarrow\prod_{n\ge0}\homzs n(\mathcal{V}(n),C^{\otimes n})\]
where $C_{\alpha}=\sum_{n\ge0}C_{\alpha_{n}}$ is countable.

The upshot is that\[
\prod_{n\ge0}\homzs n(\mathcal{V}(n),C^{\otimes n})=\dlimit\prod_{n\ge0}\homzs n(\mathcal{V}(n),C_{\alpha}^{\otimes n})\]
as $C_{\alpha}$ runs over all subcomplexes of the lattice $\mathrf{G}(C)$.\end{proof}
\begin{thm}
\label{th:coalgdirectlimit}Let $\mathcal{V}=\{\mathcal{V}(n)\}$
be an operad satisfying condition~\ref{cond:mainoperadassumption}.

If $C$ is a $\mathcal{V}$-coalgebra whose underlying chain-complex
is nearly free, then \[
C=\dlimit C_{\alpha}\]
where $\{C_{\alpha}\}$ ranges over all the countable sub-coalgebras
of $C$. \end{thm}
\begin{proof}
To prove the statement, we show that every \[
x\in C\]
is contained in a countable sub-coalgebra of $C$.

Let\[
a:C\to\prod_{n\ge0}\homzs n(\mathcal{V}(n),C^{\otimes n})\]
be the adjoint structure-map of $C$, and let $x\in C_{1}$, where
$C_{1}$ is a countable sub-chain-complex of $\forgetful C$. 

Then $a(C_{1})$ is a countable subset of $\prod_{n\ge0}\homzs n(\mathcal{V}(n),C^{\otimes n})$,
each element of which is defined by its value on the countable set
of $\zs n$-projective generators of $\{\mathcal{V}_{n}\}$ for all
$n>0$. It follows that the targets of these projective generators
are a countable set of elements\[
\{x_{j}\in C^{\otimes n}\}\]
for $n>0$. If we enumerate all of the $c_{i,j}$ in $x_{j}=c_{1,j}\otimes\cdots\otimes c_{n,j}$,
we still get a countable set. Let \[
C_{2}=C_{1}+\sum_{i,j}\ring\cdot c_{i,j}\]
This will be a countable sub-chain-complex of $\forgetful C$ that
contains $x$. By an easy induction, we can continue this process,
getting a sequence $\{C_{n}\}$ of countable sub-chain-complexes of
$\forgetful C$ with the property\[
a(C_{i})\subseteq\prod_{n\ge0}\homzs n(\mathcal{V}(n),C_{i+1}^{\otimes n})\]
arriving at a countable sub-chain-complex of $\forgetful C$\[
C_{\infty}=\bigcup_{i=1}^{\infty}C_{i}\]
that is closed under the coproduct of $C$. It is not hard to see
that the induced coproduct on $C_{\infty}$ will \emph{inherit} the
identities that make it a $\mathcal{V}$-coalgebra.\end{proof}
\begin{cor}
\label{cor:cofreedirectlimits} Let $\mathcal{V}=\{\mathcal{V}(n)\}$
be a $\Sigma$-cofibrant operad such that $\mathcal{V}(n)$ is of
finite type for all $n\ge0$. If $C$ is nearly-free, then the cofree
coalgebras \[
L_{\mathcal{V}}C,\, P_{\mathcal{V}}C,\, M_{\mathcal{V}}C,\,\mathrf{F}_{\mathcal{V}}C\]
are well-defined and\[
\left.\begin{array}{c}
L_{\mathcal{V}}C=\dlimit L_{\mathcal{V}}C_{\alpha}\\
P_{\mathcal{V}}C=\dlimit P_{\mathcal{V}}C_{\alpha}\\
M_{\mathcal{V}}C=\dlimit M_{\mathcal{V}}C_{\alpha}\\
\mathrf{F}_{\mathcal{V}}C=\dlimit\mathrf{F}_{\mathcal{V}}C_{\alpha}\end{array}\right\} \subseteq\prod_{n\ge0}\homzs n(\mathcal{V}(n),C^{\otimes n})\]
where $C_{\alpha}$ ranges over the countable sub-chain-complexes
of $C$.\end{cor}
\begin{proof}
The near-freeness of $C$ implies that the $C_{\alpha}$ are all $\integers$-free
when $\ring=\integers$, so the construction in \cite{Smith-cofree}
gives cofree coalgebras $L_{\mathcal{V}}C_{\alpha}$. 

Since (by theorem~\ref{th:coalgdirectlimit})\[
C=\dlimit C_{\alpha}\]
 where $C_{\alpha}$ ranges over countable sub-coalgebras of $C$,
we get coalgebra morphisms\[
b_{\alpha}:C_{\alpha}\to L_{\mathcal{V}}\forgetful{C_{\alpha}}\]
inducing a coalgebra morphism\[
b:C\to\dlimit L_{\mathcal{V}}\forgetful{C_{\alpha}}\]

We claim that $L_{\mathcal{V}}\forgetful C=\dlimit L_{\mathcal{V}}\forgetful{C_{\alpha}}$.
We first note that $\dlimit L_{\mathcal{V}}\forgetful{C_{\alpha}}$
depends only on $ $$\forgetful C$ and not on $C$ . If $D$ is a
$\mathcal{V}$-coalgebra with $\forgetful C=\forgetful D$ then, by
theorem~~\ref{th:coalgdirectlimit}, $D=\dlimit D_{\beta}$ where
the $D_{\beta}$ are the countable sub-coalgebras of $D$. 

We also know that, in the poset of sub-chain-complexes of $\forgetful C=\forgetful D$,
$\{\forgetful{C_{\alpha}}\}$ and $\{\forgetful{D_{\beta}}\}$ are
both cofinal. This implies the cofinality of $\{L_{\mathcal{V}}\forgetful{C_{\alpha}}\}$
and $\{L_{\mathcal{V}}\forgetful{D_{\beta}}\}$, hence\[
\dlimit L_{\mathcal{V}}\forgetful{C_{\alpha}}=\dlimit L_{\mathcal{V}}\forgetful{D_{\beta}}\]
This unique $\mathcal{V}$-coalgebra has all the categorical properties
of the cofree-coalgebra\[
L_{\mathcal{V}}\forgetful C\]
which proves the first part of the result.

The statement that\[
L_{\mathcal{V}}\forgetful C\subseteq\prod_{n\ge0}\homzs n(\mathcal{V}(n),C^{\otimes n})\]
follows from
\begin{enumerate}
\item The canonical inclusion\[
L_{\mathcal{V}}C_{\alpha}\subseteq\prod_{n\ge0}\homzs n(\mathcal{V}(n),C_{\alpha}^{\otimes n})\]
in \cite{Smith-cofree}, and
\item the fact that the hypotheses imply that\[
\prod_{n\ge0}\homzs n(\mathcal{V}(n),C^{\otimes n})=\dlimit\prod_{n\ge0}\homzs n(\mathcal{V}(n),C_{\alpha}^{\otimes n})\]
--- see proposition~\ref{prop:proddirectlimit}.
\end{enumerate}
Similar reasoning applies to $P_{\mathcal{V}}C,\, M_{\mathcal{V}}C,\,\mathrf{F}_{\mathcal{V}}C$.
\end{proof}

\subsection{Core of a module\label{sub:Core-of-a}}
\begin{lem}
\label{lem:maximalacoalg}Let $A,B\subseteq C$ be sub-coalgebras
of $C\in\mathrf{C}=\coalgcat$ or $\ircoalgcat$. Then $A+B\subseteq C$
is also a sub-coalgebra of $C$. 

In particular, given any sub-DG-module\[
M\subseteq\forgetful C\]
there exists a maximal sub-coalgebra $\core M$ --- called the core
of $M$ --- with the universal property that any sub-coalgebra $A\subseteq C$
with $\forgetful A\subseteq M$ is a sub-coalgebra of $\core M$.

This is given by\[
\alpha(\core M)=\alpha(C)\cap P_{\mathcal{V}}M\subseteq P_{\mathcal{V}}C\]
where \[
\alpha:C\to P_{\mathcal{V}}C\]
is the classifying morphism of $C$.\end{lem}
\begin{proof}
The first claim is clear --- $A+B$ is clearly closed under the coproduct
structure. This implies the second claim because we can always form
the sum of any set of sub-coalgebras contained in $M$.

The second claim follows from:

The fact that \[
\core M=\alpha^{-1}(\alpha(C)\cap P_{\mathcal{V}}M)\]
implies that it is the inverse image of a coalgebra (the intersection
of two coalgebras) under an injective map ($\alpha$), so it is a
subcoalgebra of $C$ with $\forgetful{\core M}\subseteq M$.

Given any subcoalgebra $A\subseteq C$ with $\forgetful A\subseteq M$,
the diagram \[
\xymatrix{{A}\ar@{^{(}->}[r] & {C}\ar[r]^{\alpha}\ar@{=}[rd] & {P_{\mathcal{V}}C}\ar[d]^{\epsilon}\\
{} & {} & {C}}
\]
where $\epsilon:P_{\mathcal{V}}C\to C$ is the cogeneration map, implies
that\begin{eqnarray*}
\alpha(A) & \subseteq & \alpha(C)\\
\epsilon(\alpha(A)) & \subseteq & \epsilon(P_{\mathcal{V}}M)\end{eqnarray*}
which implies that $A\subseteq\core M$, so $\core M$ has the required
universal property.\end{proof}
\begin{cor}
\label{cor:corehomotopy}Let $C\in\mathrf{C}=\coalgcat$ or $\ircoalgcat$
and $M\subseteq\forgetful C$ a sub-DG-module and suppose\[
\Phi:C\otimes I\to C\]
is a coalgebra morphism with the property that $\Phi(M\otimes I)\subseteq M$.
Then \[
\Phi(\core M\otimes I)\subseteq\core M\]
\end{cor}
\begin{proof}
The hypotheses imply that the diagrams \[
\xymatrix{{C\otimes I}\ar[r]^{\alpha\otimes1\quad}\ar@{=}[d]\ar[r] & {(L_{\mathcal{V}}C)\otimes I}\ar[d]\\
{C\otimes I}\ar[d]_{\Phi}\ar[r] & {L_{\mathcal{V}}(C\otimes I)}\ar[d]^{L_{\mathcal{V}}\Phi}\\
{C}\ar[r]_{\alpha} & {L_{\mathcal{V}}C}}
\]
and \[
\xymatrix{{(L_{\mathcal{V}}M)\otimes I}\ar[r]\ar[d] & {(L_{\mathcal{V}}C)\otimes I}\ar[d]\\
{L_{\mathcal{V}}(M\otimes I)}\ar[r]\ar[d]_{L_{\mathcal{V}}(\Phi|M\otimes I)} & {L_{\mathcal{V}}(C\otimes I}\ar[d]^{L_{\mathcal{V}}\Phi}\\
{L_{\mathcal{V}}M}\ar[r] & {L_{\mathcal{V}}C}}
\]
commute. Lemma~\ref{lem:maximalacoalg} implies the result.
\end{proof}
This allows us to construct \emph{equalizers} in categories of coalgebras
over operads:
\begin{cor}
\label{cor:coalgequalizers}If\[
f_{i}:A\to B\]
 with $i$ running over some index set, is a set of morphisms in $\mathrf{C}=\coalgcat$
or $\coalgcat$, then the equalizer of the $\{f_{i}\}$ is \[
\core M\subseteq A\]
where $M$ is the equalizer of $\forgetful{f_{i}}:\forgetful A\to\forgetful B$
in $\chaincat$.\end{cor}
\begin{rem*}
Roughly speaking, it is easy to construct coequalizers of coalgebra
morphisms and hard to construct equalizers --- since the kernel of
a morphism is not necessarily a sub-coalgebra. This is dual to what
holds for \emph{algebras} over operads.\end{rem*}
\begin{proof}
Clearly $f_{i}|\core M=f_{j}|\core M$ for all $i,j$. On the other
hand, any sub-DG-algebra with this property is contained in $\core M$
so the conclusion follows. \end{proof}
\begin{prop}
\label{pro:corecommutesintersection}Let $C\in\ircoalgcat$ and let
$\{A_{i}\}$, $i\ge0$ be a descending sequence of sub-chain-complexes
of $\forgetful C$ --- i.e., $A_{i+1}\subseteq A_{i}$ for all $i\ge0$.
Then\[
\core{\bigcap_{i=0}^{\infty}A_{i}}=\bigcap_{i=0}^{\infty}\core{A_{i}}\]
\end{prop}
\begin{proof}
Clearly, any intersection of coalgebras is a coalgebra, so\[
\bigcap_{i=0}^{\infty}\core{A_{i}}\subseteq\core{\bigcap_{i=0}^{\infty}A_{i}}\]
On the other hand\[
\forgetful{\core{\bigcap_{i=0}^{\infty}A_{i}}}\subseteq\bigcap_{i=0}^{\infty}A_{i}\subseteq A_{n}\]
for any $n>0$. Since $\core{\bigcap_{i=0}^{\infty}A_{i}}$ is a \emph{coalgebra}
whose underlying chain complex is contained in $A_{n}$, we must actually
have\[
\forgetful{\core{\bigcap_{i=0}^{\infty}A_{i}}}\subseteq\forgetful{\core{A_{n}}}\]
which implies that\[
\core{\bigcap_{i=0}^{\infty}A_{i}}\subseteq\bigcap_{i=0}^{\infty}\core{A_{i}}\]
and the conclusion follows.\end{proof}
\begin{defn}
\label{def:aveeb}Let $A$ and $B$ be objects of $\mathrf{C}=\coalgcat$
or $\ircoalgcat$ and define $A\vee B$ to be the \emph{push out}
in the diagram \[
\xymatrix{{\bullet}\ar[r]\ar[d] & {B}\ar@{.>}[d]\\
{A}\ar@{.>}[r] & {A\vee B}}
\]
 where where $\bullet$ denotes the initial object in $\mathrf{C}$
--- see definition~\ref{def:pointedirredcat}. 
\end{defn}

\subsection{Categorical products\label{sub:Categorical-products}}

We can use cofree coalgebras to explicitly construct the categorical
product in $\ircoalgcat$ or $\coalgcat$:
\begin{defn}
\label{def:catprod}Let $A_{i}$, $i=0,1$ be objects of $\mathrf{C}=\coalgcat$
or $\ircoalgcat$. Then \[
A_{0}\boxtimes A_{1}=\core{M_{0}\cap M_{1}}\subseteq Z=\begin{cases}
L_{\mathcal{V}}(\forgetful{A_{0}}\oplus\forgetful{A_{1}}) & \mathrm{if\,}\mathrf{C}=\coalgcat\\
\pcoalg V{\forgetful{A_{0}}\oplus\forgetful{A_{1}}} & \mathrm{if\,}\mathrf{C}=\ircoalgcat\end{cases}\]
where\[
M_{i}=p_{i}^{-1}(\forgetful{\im A_{i}})\]
under the projections\[
p_{i}:Z\to\begin{cases}
L_{\mathcal{V}}\forgetful{A_{i}} & \mathrm{if\,}\mathrf{C}=\coalgcat\\
\pcoalgf V{A_{i}} & \mathrm{if\,}\mathrf{C}=\ircoalgcat\end{cases}\]
 induced by the canonical maps $\forgetful{A_{0}}\oplus\forgetful{A_{1}}\to\forgetful{A_{i}}$.
The $\im A_{i}$ are images under the canonical morphisms\[
A_{i}\to\left\{ \begin{array}{cc}
L_{\mathcal{V}}\forgetful{A_{i}} & \mathrm{if\,}\mathrf{C}=\coalgcat\\
\pcoalgf V{A_{i}} & \mathrm{if\,}\mathrf{C}=\ircoalgcat\end{array}\right\} \to Z\]
classifying coalgebra structures --- see definition~\ref{def:cofreecoalgebra}.\end{defn}
\begin{rem*}
By identifying the $A_{i}$ with their canonical images in $Z$, we
get canonical projections to the factors\[
A_{0}\boxtimes A_{1}\to A_{i}\]
\end{rem*}
\begin{prop}
\label{pro:stricthomotopyproducts}Let $F:A\otimes I\to A$ and $G:B\otimes I\to B$
be strict homotopies in $\mathrf{C}=\ircoalgcat$ or $\coalgcat$
(see definition~\ref{def:strictequivalence}). Then there is a strict
homotopy\[
(A\boxtimes B)\otimes I\xrightarrow{F\hat{\boxtimes}G}A\boxtimes B\]
 that makes the diagrams \begin{equation}
\xymatrix{{(A\boxtimes B)\otimes I}\ar[d]\ar[r]^{\quad F\hat{\boxtimes}G} & {A\boxtimes B}\ar[d]\\
{A\otimes I}\ar[r]_{F} & {A}}
\label{eq:boxtimesdia1}\end{equation}
 and \begin{equation}
\xymatrix{{(A\boxtimes B)\otimes I}\ar[d]\ar[r]^{\quad F\hat{\boxtimes}G} & {A\boxtimes B}\ar[d]\\
{B\otimes I}\ar[r]_{G} & {B}}
\label{eq:boxtimesdia2}\end{equation}
commute. If\begin{eqnarray*}
f_{1},f_{2}:A & \to & A'\\
g_{1},g_{2}:B & \to & B'\end{eqnarray*}
 are strictly homotopic morphisms with respective strict homotopies\begin{eqnarray*}
F:A\otimes I & \to & A'\\
G:B\otimes I & \to & B'\end{eqnarray*}
then $F\hat{\boxtimes}G$ is a strict homotopy between $f_{1}\boxtimes g_{1}$
and $f_{2}\boxtimes g_{2}$.

Consequently, if $f:A\to A'$ and $g:B\to B'$ are strict equivalences,
then \[
f\boxtimes g:A\boxtimes B\to A'\boxtimes B'\]
 is a strict equivalence.\end{prop}
\begin{proof}
The projections\[
\xymatrix{{A\boxtimes B}\ar[d]\ar[r] & {B}\\
{A} & {}}
\]
induce projections\[
\xymatrix{{(A\boxtimes B)\otimes I}\ar[d]\ar[r] & {B\otimes I}\\
{A\otimes I} & {}}
\]
and the composite\[
(A\boxtimes B)\otimes I\to(A\otimes I)\boxtimes(B\otimes I)\xrightarrow{F\boxtimes G}A\boxtimes B\]
 satisfies the first part of the statement. 

Note that diagrams~\ref{eq:boxtimesdia1} and \ref{eq:boxtimesdia2}
--- and the fact that maps to $A'\boxtimes B'$ are \emph{uniquely
determined} by their composites with the projections $A'\boxtimes B'\rightrightarrows A',B'$
(the defining universal property of $\boxtimes$) --- implies that
$F\hat{\boxtimes}G$ is a strict homotopy between $f_{1}\boxtimes g_{1}$
and $f_{2}\boxtimes g_{2}$. The final statement is also clear.
\end{proof}
In like fashion, we can define categorical \emph{fibered} products
of coalgebras:
\begin{defn}
\label{def:catfiberedproduct}Let \[
\xymatrix{{F}\ar[r]\ar[d]_{i} & {A}\ar[d]^{f}\\
{B}\ar[r] & {C}}
\]
be a diagram in $\coalgcat$ or $\ircoalgcat$. Then the \emph{fibered
product with respect to this diagram, $\mathit{A}\boxtimes^{\mathit{C}}\mathit{B}$,}
is defined to be the equalizer\[
F\to A\boxtimes B\rightrightarrows C\]
by the maps induced by the projections $A\boxtimes B\to A$ and $A\boxtimes B\to B$
composed with the maps in the diagram.
\end{defn}
We have an analogue to proposition~\ref{pro:stricthomotopyproducts}:
\begin{prop}
\label{pro:strictfiberedprodinvar}Let $A\xrightarrow{f}B\xleftarrow{g}C$,
$A'\xrightarrow{f'}B'\xleftarrow{g'}C'$ be diagrams in $\mathrf{C}=\ircoalgcat$
or $\coalgcat$ and let \[
\xymatrix{{A\otimes I}\ar[r]^{f\otimes1}\ar[d]_{H_{A}} & {B\otimes I}\ar[d]^{H_{B}} & {C\otimes I}\ar[l]_{g\otimes1}\ar[d]^{H_{C}}\\
{A'}\ar[r]_{f'} & {B'} & {C}\ar[l]^{g'}}
\]
commute, where the $H_{\alpha}$ are strict homotopies. Then there
exists a strict homotopy\[
(A\boxtimes^{B}C)\otimes I\xrightarrow{H_{A}\hat{\boxtimes}^{H_{B}}H_{C}}A'\boxtimes^{B'}C'\]
between the morphisms\[
(H_{A}|A\otimes p_{i})\boxtimes(H_{C}\otimes p_{i}):A\boxtimes^{B}C\to A'\boxtimes^{B'}C'\]
for $i=0,1$.\end{prop}
\begin{proof}
The morphism $H_{A}\hat{\boxtimes}^{H_{B}}H_{C}$ is constructed exactly
as in proposition~\ref{pro:stricthomotopyproducts}. The conclusion
follows by the same reasoning used to prove the final statement of
that result.\end{proof}
\begin{prop}
\label{prop:freepullback}Let $U,V$ and $W$ be objects of $\chaincat$
and let $Z$ be the fibered product of \[
\xymatrix{{} & {V}\ar[d]^{g}\\
{U}\ar[r]_{f} & {W}}
\]
 in $\chaincat$ --- i.e., $W$ is the equalizer \[
Z\to U\oplus V\rightrightarrows W\]
in $\chaincat$. Then $P_{\mathcal{V}}Z$ is the fibered product of
\begin{equation}
\xymatrix{{} & {P_{\mathcal{V}}V}\ar[d]^{P_{\mathcal{V}}g}\\
{P_{\mathcal{V}}U}\ar[r]_{P_{\mathcal{V}}f} & {P_{\mathcal{V}}W}}
\label{dia:cofreefibered}\end{equation}
 in $\ircoalgcat$ and $L_{\mathcal{V}}Z$ is the fibered product
of \begin{equation}
\xymatrix{{} & {L_{\mathcal{V}}V}\ar[d]^{L_{\mathcal{V}}g}\\
{L_{\mathcal{V}}U}\ar[r]_{L_{\mathcal{V}}f} & {L_{\mathcal{V}}W}}
\end{equation}
 in $\coalgcat$.\end{prop}
\begin{proof}
We prove this in the pointed irreducible case. The other case follows
by an analogous argument.

The universal properties of cofree coalgebras imply that $P_{\mathcal{V}}(U\oplus V)=P_{\mathcal{V}}U\boxtimes P_{\mathcal{V}}V$.
Suppose $F$ is the fibered product of diagram~\ref{dia:cofreefibered}.
Then\[
P_{\mathcal{V}}Z\subseteq F\]
On the other hand, the composite\[
F\to P_{\mathcal{V}}U\boxtimes P_{\mathcal{V}}V=P_{\mathcal{V}}(U\oplus V)\to U\oplus V\]
 where the rightmost map is the co-generating map, has composites
with $f$ and $g$ that are equal to each other --- so it lies in
$Z\subseteq U\oplus V$. This induces a \emph{unique} coalgebra morphism\[
j:F\to P_{\mathcal{V}}Z\]
 left-inverse to the inclusion \[
i:P_{\mathcal{V}}Z\subseteq F\]
The uniqueness of induced maps to cofree coalgebras implies that $j\circ i=i\circ j=1$.
\end{proof}

\subsection{Limits and colimits\label{sub:Limits-and-colimits}}

We can use cofree coalgebras and adjointness to the forgetful functors
to define categorical limits and colimits in $\ircoalgcat$ and $\coalgcat$.

Categorical reasoning implies that
\begin{prop}
\label{pro:ilimitfibereddia}Let \[
\xymatrix{{} & {\{B_{i}\}}\ar[d]^{b_{i}}\\
{\{A_{i}\}}\ar[r]_{a_{i}} & {\{C_{i}\}}}
\]
be a diagram in $\iircoalgcat$ or $\icoalgcat$. Then\[
\ilimit(A_{i}\boxtimes^{C_{i}}B_{i})=(\ilimit A_{i})\boxtimes^{(\ilimit C_{i})}(\ilimit B_{i})\]
 See definition~\ref{def:catfiberedproduct} for the fibered product
notation.
\end{prop}
Theorem~\ref{th:cofreelimits} implies that colimits in $\ircoalgcat$
or $\coalgcat$ are the same as colimits of underlying chain-complexes.
The corresponding statement for limits is not true except in a special
case:
\begin{prop}
\label{pro:ilimitintersection}Let $\{C_{i}\}\in\iircoalgcat$ or
$\icoalgcat$ and suppose that all of its morphisms are injective.
Then\[
\forgetful{\ilimit C_{i}}=\ilimit\forgetful{C_{i}}\]
\end{prop}
\begin{rem*}
In this case, the limit is an intersection of coalgebras. This result
says that to get the limit of $\{C_{i}\}$, one
\begin{enumerate}
\item forms the limit of the underlying chain-complexes (i.e., the intersection)
and 
\item equips that with the coalgebra structure in induced by its inclusion
into any of the $C_{i}$
\end{enumerate}
That this constructs the limit follows from the \emph{uniqueness}
of limits.\end{rem*}
\begin{defn}
\label{def:normalizedinversesystem}Let $A=\{A_{i}\}\in\iircoalgcat$.
Then define the \emph{normalization of} $A$, denoted $\hat{A}=\{\hat{A}_{i}\}$,
as follows:
\begin{enumerate}
\item Let $V=\pcoalg V{\ilimit\forgetful{A_{i}}}$ with canonical maps \[
q_{n}:\pcoalg V{\ilimit\forgetful{A_{i}}}\to\pcoalgf V{A_{n}}\]
for all $n>0$.
\item Let $f_{n}:A_{n}\to P_{\mathcal{V}}(\forgetful{A_{n}})$ be the coalgebra
classifying map --- see definition~\ref{def:cofreecoalgebra}.
\end{enumerate}
Then $\hat{A}_{n}=\core{q_{n}^{-1}(f_{n}(A_{n}))}$, and $\hat{A}_{n+1}\subseteq\hat{A}_{n}$
for all $n>0$. Define $\hat{A}=\{\hat{A}_{n}\}$, with the injective
structure maps defined by inclusion.

If $A=\{A_{i}\}\in\icoalgcat$ then the corresponding construction
holds, where we consistently replace $P_{\mathcal{V}}(*)$ by $L_{\mathcal{V}}(*)$.
\end{defn}
Normalization reduces the \emph{general} case to the case dealt with
in proposition~\ref{pro:ilimitintersection}.
\begin{cor}
\label{cor:generalcoalglimits}Let $C=\{g_{i}:C_{i}\to C_{i-1}\}$
in $\iircoalgcat$ or $\icoalgcat$. Then\[
\ilimit C_{i}=\ilimit\hat{C}_{i}\]
where $\{\hat{C}_{i}\}$ is the normalization of $\{C_{i}\}$. In
particular, if $C$ is in $\ircoalgcat$ \[
\ilimit C_{i}=\core{\bigcap_{i=0}^{\infty}\forgetful{p_{i}}^{-1}\forgetful{\alpha_{i}}(\forgetful{C_{i}})}\subseteq P_{\mathcal{V}}(\ilimit\forgetful{C_{i}})\]
where $p_{i}:P_{\mathcal{V}}(\ilimit\forgetful{C_{i}})\to P_{\mathcal{V}}(\forgetful{C_{i}})$
and $\alpha_{n}:C_{n}\to P_{\mathcal{V}}(\forgetful{C_{i}})$ are
as in definition~\ref{def:normalizedinversesystem}, and the corresponding
statement holds if $C$ is in $\icoalgcat$ with $P_{\mathcal{V}}(*)$
replaced by $L_{\mathcal{V}}(*)$ .\end{cor}
\begin{proof}
Assume the notation of definition~\ref{def:normalizedinversesystem}.
Let \[
f_{i}:C_{i}\to\left\{ \begin{array}{c}
\pcoalgf V{C_{i}}\\
L_{\mathcal{V}}\forgetful{C_{i}}\end{array}\right\} \]
be the classifying maps in $\ircoalgcat$ or $\coalgcat$, respectively
--- see definition~\ref{def:cofreecoalgebra}. We deal with the case
of the category $\ircoalgcat$ --- the other case is entirely analogous.
Let\[
q_{n}:\pcoalg V{\ilimit\forgetful{C_{i}}}\to\pcoalgf V{C_{n}}\]
be induced by the canonical maps $\ilimit\forgetful{C_{i}}\to\forgetful{C_{n}}$.

We verify that \[
X=\core{\bigcap_{i=0}^{\infty}\forgetful{q_{i}}^{-1}\forgetful{f_{i}}(\forgetful{C_{i}})}=\ilimit\hat{C}_{i}\]
has the category-theoretic properties of an inverse limit. We must
have morphisms\[
p_{i}:X\to C_{i}\]
 making the diagrams\begin{equation}
\xymatrix{{X}\ar[r]^{p_{i}}\ar[rd]_{p_{i-1}} & {C_{i}}\ar[d]^{g_{i}}\\
{} & {C_{i-1}}}
\label{eq:inversediagram}\end{equation}
 commute for all $i>0$. Define $p_{i}=f_{i}^{-1}\circ q_{i}:X\to C_{i}$
--- using the fact that the classifying maps $f_{i}:C_{i}\to P_{\mathcal{V}}\forgetful{C_{i}}$
are always injective (see \cite{Smith-cofree} and the definition~\ref{def:cofreecoalgebra}).
The commutative diagrams \[
\xymatrix{{C_{i}}\ar[d]_{g_{i}}\ar[r]^{\alpha_{i}} & {P_{\mathcal{V}}\forgetful{C_{i}}}\ar[d]^{P_{\mathcal{V}}\forgetful{g_{i}}}\\
{C_{i-1}}\ar[r]_{\alpha_{i-1}} & {P_{\mathcal{V}}\forgetful{C_{i-1}}}}
\]
 and \[
\xymatrix{{\ilimit P_{\mathcal{V}}\forgetful{C_{i}}}\ar[r]^{p_{i}}\ar[rd]_{p_{i-1}} & {P_{\mathcal{V}}\forgetful{C_{i}}}\ar[d]^{P_{\mathcal{V}}\forgetful{g_{i}}}\\
{} & {P_{\mathcal{V}}\forgetful{C_{i-1}}}}
\]
together imply the commutativity of the diagram with the diagrams
\ref{eq:inversediagram}. Consequently, $X$ is a \emph{candidate}
for being the inverse limit, $\ilimit C_{i}$.

We must show that any other candidate $Y$ possesses a \emph{unique}
morphism $Y\to X$, making appropriate diagrams commute. Let $Y$
be such a candidate. The morphism of inverse systems defined by classifying
maps (see definition~\ref{def:cofreecoalgebra})\[
C_{i}\to\pcoalgf V{C_{i}}\]
 implies the existence of a \emph{unique} morphism \[
Y\to\ilimit P_{\mathcal{V}}\forgetful{C_{i}}=P_{\mathcal{V}}\ilimit\forgetful{C_{i}}\]
The commutativity of the diagrams\[
\xymatrix{{Y}\ar[r]\ar[d] & {P_{\mathcal{V}}\ilimit\forgetful{C_{i}}}\ar[d]^{p_{i}}\\
{C_{i}}\ar[r]_{\alpha_{i}} & {P_{\mathcal{V}}\forgetful{C_{i}}}}
\]
 for all $i\ge0$ implies that $\forgetful{\im Y}\subseteq\forgetful{p_{i}}^{-1}\forgetful{\alpha_{i}}(\forgetful{C_{i}})$.
Consequently \[
\forgetful{\im Y}\subseteq\bigcap_{i=0}^{\infty}\forgetful{p_{i}}^{-1}\forgetful{\alpha_{i}}(\forgetful{C_{i}})\]
Since $Y$ is a coalgebra, its image must lie within the maximal sub-coalgebra
contained within $\bigcap_{i=0}^{\infty}\forgetful{p_{i}}^{-1}\forgetful{\alpha_{i}}(\forgetful{C_{i}})$,
namely $X=\core{\bigcap_{i=0}^{\infty}\forgetful{p_{i}}^{-1}\forgetful{\alpha_{i}}(\forgetful{C_{i}})}$.
This proves the first claim. Proposition~\ref{pro:corecommutesintersection}
implies that $X=\bigcap_{i=0}^{\infty}\hat{C}_{i}=\ilimit\hat{C}_{i}$.\end{proof}
\begin{lem}
\label{lem:tensorinvliminjective}Let $\{g_{i}:C_{i}\to C_{i-1}\}$
be an inverse system in $\chaincat$. If $n>0$ is an integer, then
the natural map\[
\left(\ilimit C_{i}\right)^{\otimes n}\to\ilimit C_{i}^{\otimes n}\]
is injective.\end{lem}
\begin{proof}
Let $A=\ilimit C_{i}$ and $p_{i}:A\to C_{i}$ be the natural projections.
If \[
W_{k}=\ker p_{k}^{\otimes n}:\left(\ilimit C_{i}\right)^{\otimes n}\to C_{k}^{\otimes n}\]
  we will show that\[
\bigcap_{k=1}^{\infty}W_{k}=0\]
If $K_{i}=\ker p_{i}$, then\[
\bigcap_{i=1}^{\infty}K_{i}=0\]
and \[
W_{i}=\sum_{j=1}^{n}\underbrace{A\otimes\cdots\otimes K_{i}\otimes\cdots\otimes A}_{j^{\text{th}}\,\text{position}}\]
 Since all modules are nearly-free, hence, flat (see remark~\ref{rem:nearlyfreeimpliesflat}),
we have\[
W_{k+1}\subseteq W_{k}\]
 for all $k$, and \[
\bigcap_{i=1}^{m}W_{i}=\sum_{j=1}^{n}\underbrace{A\otimes\cdots\otimes\left(\bigcap_{i=1}^{m}K_{i}\right)\otimes\cdots\otimes A}_{j^{\text{th}}\,\text{position}}\]
from which the conclusion follows.\end{proof}
\begin{prop}
\label{pro:injectioninvercofree}Let $\{C_{i}\}\in\iircoalgcat$,
and suppose $\mathcal{V}=\{\mathcal{V}(n)\}$ is a $\Sigma$-cofibrant
operad with $\mathcal{V}(n)$ of finite type for all $n\ge0$. Then
the projections \[
\forgetful{\pcoalg V{\ilimit\forgetful{C_{i}}}}\to\forgetful{\pcoalgf V{C_{n}}}\]
 for all $n>0,$ induce a canonical injection \[
\mu:\forgetful{\pcoalg V{\ilimit\forgetful{C_{i}}}}\hookrightarrow\ilimit\forgetful{\pcoalgf V{C_{i}}}\]
In addition, the fact that the structure maps \[
\alpha_{i}:C_{i}\to P_{\mathcal{V}}(\forgetful{C_{i}})\]
of the $\{C_{i}\}$ are coalgebra morphisms implies the existence
of an injective $\chaincat$-morphism\[
\hat{\alpha}:\ilimit\forgetful{C_{i}}\hookrightarrow\ilimit\forgetful{\pcoalgf V{C_{i}}}\]

Corresponding statements hold for $\icoalgcat$ and the functors $L_{\mathcal{V}}(*)$.\end{prop}
\begin{proof}
We must prove that \[
\mu:\forgetful{\pcoalg V{\ilimit\forgetful{C_{i}}}}\to\ilimit\forgetful{\pcoalgf V{C_{i}}}\]
is injective. Let $K=\ker\mu$. Then \[
K\subset\forgetful{\pcoalg V{\ilimit\forgetful{C_{i}}}}\subseteq\prod_{n\ge0}\homzs n(\mathcal{V}(n),D^{\otimes n})\]
where $D=\ilimit\forgetful{C_{i}}$ (see \cite{Smith-cofree}). If
$n\ge0$, let\[
p_{n}:\prod_{n\ge0}\homzs n(\mathcal{V}(n),D^{\otimes n})\to\homzs n(\mathcal{V}(n),D^{\otimes n})\]
denote the canonical projections. The diagrams \[
\xyC{80pt}\xymatrix{\prod_{n\ge0}\homzs n(\mathcal{V}(n),D^{\otimes n})\ar[r]^{{\displaystyle \prod\homz(1,b_{k}^{\otimes n})}}\ar[d]_{p_{n}} & {\displaystyle \prod_{n\ge0}\homzs n(\mathcal{V}(n),C_{k}^{\otimes n})}\ar[d]^{q_{n}}\\
{\homzs n(\mathcal{V}(n),D^{\otimes n})}\ar[r]_{\homz(1,b_{k}^{\otimes n})} & {\homzs n(\mathcal{V}(n),C_{k}^{\otimes n})}}
\]
 commute for all $k$ and $n\ge0$, where $q_{n}$ is the counterpart
of $p_{n}$ and $b_{k}:\ilimit\forgetful{C_{i}}\to\forgetful{C_{k}}$
is the canonical map. It follows that \[
p_{k}(K)\subseteq\ker\homz(1,b_{k}^{\otimes n})\]
for all $n\ge0$$ $, or\[
p_{k}(K)\subseteq\bigcap_{k>0}\ker\homz(1,b_{k}^{\otimes n})\]
 We claim that\[
\bigcap_{n>0}\ker\homz(1,b_{k}^{\otimes n})=\homz(1,\bigcap_{k>0}\ker b_{k}^{\otimes n})\]

The equality on the left follows from the left-exactness of $\homz$
and filtered limits (of chain-complexes). The equality on the right
follows from the fact that
\begin{enumerate}
\item $\bigcap_{k>0}\ker b_{k}=0$
\item the left exactness of $\otimes$ for $\ring$-\emph{flat} modules
(see remark~\ref{rem:nearlyfreeimpliesflat}).
\item Lemma~\ref{lem:tensorinvliminjective}.
\end{enumerate}
It follows that $p_{n}(K)=0$ for all $n\ge0$ and $K=0$.

The map\[
\hat{\alpha}:\ilimit\forgetful{C_{i}}\hookrightarrow\ilimit\forgetful{\pcoalgf V{C_{i}}}\]
is induced by \emph{classifying maps} of the coalgebras $\{C_{i}\}$,
which induce a morphism of limits because the structure maps $C_{k}\to C_{k-1}$
are coalgebra morphisms, making the diagrams \[
\xymatrix{{C_{k}}\ar[r]\ar[d] & {C_{k-1}}\ar[d]\\
{\pcoalgf V{C_{k}}}\ar[r] & {\pcoalgf V{C_{k-1}}}}
\]
commute for all $k>0$.\end{proof}
\begin{cor}
\label{cor:intersectionalreadycoalg}Let $C=\{g_{i}:C_{i}\to C_{i-1}\}\in\iircoalgcat$,
and suppose $\mathcal{V}=\{\mathcal{V}(n)\}$ is a $\Sigma$-cofibrant
operad with $\mathcal{V}(n)$ of finite type for all $n\ge0$. Then\begin{equation}
\ilimit C_{i}=\bigcap_{i=1}^{\infty}q_{i}^{-1}(\alpha_{i}(C_{i}))\subseteq L_{\mathcal{V}}(\ilimit\forgetful{C_{i}})\label{eq:bigintersectioncoalg}\end{equation}
with the coproduct induced from $L_{\mathcal{V}}(\ilimit\forgetful{C_{i}})$,
and where \[
q_{i}:L_{\mathcal{V}}(\ilimit\forgetful{C_{i}})\to L_{\mathcal{V}}(\forgetful{C_{i}})\]
 is the projection and \begin{equation}
\alpha_{i}:C_{i}\to L_{\mathcal{V}}(\forgetful{C_{i}})\label{eq:alphai}\end{equation}
is the classifying map, for all $i$. In addition, the sequence\begin{multline}
0\to\forgetful{\ilimit C_{i}}\to\ilimit\forgetful{C_{i}}\xrightarrow{\hat{\alpha}}\frac{\ilimit\forgetful{L_{\mathcal{V}}(\forgetful{C_{i}})}}{\mu(\forgetful{L_{\mathcal{V}}(\ilimit\forgetful{C_{i}})})}\\
\to\frac{\ilimit(\forgetful{L_{\mathcal{V}}(\forgetful{C_{i}})/\forgetful{\alpha_{i}(C_{i})}})}{\im\mu(\forgetful{L_{\mathcal{V}}(\ilimit\forgetful{C_{i}})})}\to\ilimit^{1}\forgetful{C_{i}}\to0\label{eq:sixterminvlimexact}\end{multline}
is exact in $\chaincat$, where the injection\[
\forgetful{\ilimit C_{i}}\to\ilimit\forgetful{C_{i}}\]
is induced by the projections\[
p_{i}:\ilimit C_{i}\to C_{i}\]

and \[
\hat{\alpha}:\ilimit\forgetful{C_{i}}\to\ilimit\forgetful{L_{\mathcal{V}}(\forgetful{C_{i}})}\]
 is induced by the $\{\alpha_{i}\}$ in equation~\ref{eq:alphai}.
The map\[
\mu:\forgetful{L_{\mathcal{V}}(\ilimit\forgetful{C_{i}}}\hookrightarrow\ilimit\forgetful{L_{\mathcal{V}}(\forgetful{C_{i}}}\]
is constructed in Proposition~\ref{pro:injectioninvercofree}.

If $C\in\icoalgcat$, then the corresponding statements apply, where
$L_{\mathcal{V}}(*)$ is replaced by $P_{\mathcal{V}}(*)$.\end{cor}
\begin{rem*}
The first statement implies that the use of the $\core *$-functor
in corollary~\ref{cor:generalcoalglimits} is unnecessary --- at
least if $\mathcal{V}$ is projective in the sense defined above.

The remaining statements imply that $\ilimit C_{i}$ is the \emph{largest}
sub-chain-complex of $\ilimit\forgetful{C_{i}}$ upon which one can
define a coproduct that is compatible with the maps\[
\ilimit C_{i}\to C_{i}\]
\end{rem*}
\begin{proof}
First, consider the projections\[
q_{i}:L_{\mathcal{V}}(\ilimit\forgetful{C_{i}})\to L_{\mathcal{V}}(\forgetful{C_{i}})\]
The commutativity of the diagram \[
\xymatrix{{L_{\mathcal{V}}(\ilimit\forgetful{C_{i}})}\ar@{^{(}->}[r]\ar[d]_{q_{i}} & {\ilimit L_{\mathcal{V}}(C_{i})}\ar[d]\\
{L_{\mathcal{V}}(C_{i})}\ar@{=}[r] & {L_{\mathcal{V}}(C_{i})}}
\]
 implies that\[
\ilimit\ker q_{i}=\bigcap_{i=1}^{\infty}\ker q_{i}=0\]
Now, consider the exact sequence\[
0\to\ker q_{i}\to q_{i}^{-1}(\alpha_{i}(C_{i}))\to\forgetful{C_{i}}\to0\]
and pass to inverse limits. We get the standard 6-term exact sequence
for inverse limits (of $\integers$-modules):\begin{multline}
0\to\ilimit\ker q_{i}\to\ilimit q_{i}^{-1}(\alpha_{i}(\forgetful{C_{i}}))\to\ilimit\forgetful{C_{i}}\\
\to\ilimit^{1}\ker q_{i}\to\ilimit^{1}q_{i}^{-1}(\alpha_{i}(C_{i}))\to\ilimit^{1}\forgetful{C_{i}}\to0\label{eq:sixtermexactseq}\end{multline}
which, with the fact that $\ilimit\ker q_{i}=0$, implies that\[
\bigcap_{i=1}^{\infty}q_{i}^{-1}(\alpha_{i}(C_{i}))=\ilimit q_{i}^{-1}(\alpha_{i}(C_{i}))\hookrightarrow\ilimit\forgetful{C_{i}}\]
The conclusion follows from the fact that\[
\ilimit C_{i}=\core{\bigcap_{i=1}^{\infty}q_{i}^{-1}(\alpha_{i}(C_{i}))}\subseteq\bigcap_{i=1}^{\infty}q_{i}^{-1}(\alpha_{i}(C_{i}))\]

It remains to prove the claim in equation~\ref{eq:bigintersectioncoalg},
which amounts to showing that\[
J=\bigcap_{i=1}^{\infty}q_{i}^{-1}(\alpha_{i}(C_{i}))\subseteq L_{\mathcal{V}}(\ilimit\forgetful{C_{i}})\]
is closed under the coproduct of $L_{\mathcal{V}}(\ilimit\forgetful{C_{i}})$
--- i.e., it is a coalgebra even \emph{without} applying the $\core *$-functor.
If $n\ge0$, consider the diagram \begingroup\xymatrixcolsep{7pt}\small

\[
\xymatrix{{q_{j}^{-1}(\alpha_{j}(C_{j}))}\ar@{^{(}->}[d]\ar[r]^{\alpha_{j}^{-1}\circ q_{j}} & C_{j}\ar@{^{(}->}[d]_{\alpha_{j}}\ar@/^{6pc}/[dddd]^{c_{n,j}}\\
{L_{\mathcal{V}}(\ilimit\forgetful{C_{i}})}\ar[d]_{\hat{\delta}_{n}}\ar[r]_{q_{j}} & {L_{\mathcal{V}}(\forgetful{C_{j}})}\ar[d]_{\delta_{j,n}}\\
{\homzs n(\mathcal{V}(n),L_{\mathcal{V}}(\ilimit\forgetful{C_{i}})^{\otimes n}))}\ar@{^{(}->}[d]_{\hat{\mu}_{n}}\ar[r]^{\,\, r_{n,j}} & {\homzs n(\mathcal{V}(n),(L_{\mathcal{V}}(\forgetful{C_{j}})^{\otimes n}))}\\
\homzs n(\mathcal{V}(n),(\ilimit L_{\mathcal{V}}(\forgetful{C_{i}}))^{\otimes n})\ar@{^{(}->}[d]\ar[ru]_{\quad\homz(1,p_{j}^{\otimes n})}\\
\prod_{n\ge0}\homzs n(\mathcal{V}(n),L_{\mathcal{V}}(\forgetful{C_{i}})^{\otimes n})\ar@/_{2pc}/[uur]_{\pi_{j}} & \homzs n(\mathcal{V}(n),C_{j}^{\otimes n})\ar[uu]_{s_{n,j}}}
\]
 \endgroup where:
\begin{enumerate}
\item the $\delta_{i}$ and $\hat{\delta}$-maps are coproducts and the
$\alpha_{i}$ are coalgebra morphisms. 
\item $r_{n,j}=\homz(1,q_{j}^{\otimes n})$,
\item The map $\hat{\mu}_{n}$ is defined by\begin{multline*}
\hat{\mu}_{n}=\homz(1,\mu^{\otimes n}):\homzs n(\mathcal{V}(n),(L_{\mathcal{V}}(\ilimit\forgetful{C_{i}})^{\otimes n}))\\
\hookrightarrow\homzs n(\mathcal{V}(n),(\ilimit L_{\mathcal{V}}(\forgetful{C_{i}})^{\otimes n})\end{multline*}
 where $\mu:L_{\mathcal{V}}(\ilimit\forgetful{C_{i}})\hookrightarrow\ilimit L_{\mathcal{V}}(\forgetful{C_{i}})$
is the map defined in Proposition~\ref{pro:injectioninvercofree}.
\item $s_{n,j}=\homz(1,\alpha_{j}^{\otimes n})$, and $\alpha_{j}:C_{j}\to L_{\mathcal{V}}(\forgetful{C_{j}})$
is the classifying map.
\item $p_{j}:\ilimit L_{\mathcal{V}}(\forgetful{C_{i}})\to L_{\mathcal{V}}(\forgetful{C_{j}})$
is the canonical projection.
\item $c_{n,j}:C_{j}\to\homzs n(\mathcal{V}(n),C_{j}^{\otimes n})$ is the
coproduct.
\end{enumerate}
This diagram and the projectivity of $\{\mathcal{V}(n)_{*}\}$ and
the near-freeness of $L_{\mathcal{V}}(\ilimit\forgetful{C_{i}})$
(and flatness: see remark~\ref{rem:nearlyfreeimpliesflat}) implies
that \[
\hat{\delta}_{n}\left(q_{j}^{-1}(\alpha_{j}(C_{j}))\right)\subseteq\homzs n(\mathcal{V}(n),L_{n,j})\]
where $L_{n,j}=q_{j}^{-1}(\alpha_{j}(C_{j}))^{\otimes n}+\ker r_{n,j}$
and \[
\bigcap_{j=1}^{\infty}L_{n,j}=\left(\bigcap_{j=1}^{\infty}q_{j}^{-1}(\alpha_{j}(C_{j}))\right)^{\otimes n}+\ker\hat{\mu}_{n}=J^{\otimes n}\]
so $J$ is closed under the coproduct for $L_{\mathcal{V}}(\ilimit\forgetful{C_{i}})$. 

Now, we claim that the exact sequence~\ref{eq:sixterminvlimexact}
is just \ref{eq:sixtermexactseq} in another form --- we have expressed
the $\ilimit^{1}$ terms as quotients of limits of \emph{other} terms.

The exact sequences\[
0\to\ker q_{k}\to\forgetful{L_{\mathcal{V}}(\ilimit\forgetful{C_{i}})}\xrightarrow{q_{k}}\forgetful{L_{\mathcal{V}}(\forgetful{C_{k}})}\to0\]
for all $k$, induces the sequence of limits \begingroup\xymatrixcolsep{8pt} 
\[
\xymatrix{{0}\ar[r] & {\ilimit\ker q_{k}}\ar[r]\ar@{=}[d] & {\forgetful{L_{\mathcal{V}}(\ilimit\forgetful{C_{i}})}}\ar[r] & {\ilimit\forgetful{L_{\mathcal{V}}(\forgetful{C_{i}})}}\ar[r] & {\ilimit^{1}\ker q_{k}}\ar[r] & {0}\\
 & {0}}
\]
\endgroup which implies that\[
\ilimit^{1}\ker q_{k}=\frac{\ilimit\forgetful{L_{\mathcal{V}}(\forgetful{C_{i}})}}{\forgetful{L_{\mathcal{V}}(\ilimit\forgetful{C_{i}})}}\]
In like fashion, the exact sequences\[
0\to q_{k}^{-1}(\alpha_{k}(C_{k}))\to\forgetful{L_{\mathcal{V}}(\ilimit\forgetful{C_{i}})}\to\forgetful{L_{\mathcal{V}}(\forgetful{C_{k}})/\alpha_{k}(C_{k})}\to0\]
 imply that\[
\ilimit^{1}q_{i}^{-1}(\alpha_{i}(C_{i}))=\frac{\ilimit(\forgetful{L_{\mathcal{V}}(\forgetful{C_{i}})/\forgetful{\alpha_{i}(C_{i})}})}{\im\forgetful{L_{\mathcal{V}}(\ilimit\forgetful{C_{i}})}}\]
\end{proof}
\begin{cor}
\label{cor:inverselimitintersectiontwo}Let $\{C_{i}\}\in\iircoalgcat$,
and suppose $\mathcal{V}=\{\mathcal{V}(n)\}$ is a $\Sigma$-cofibrant
operad with $\mathcal{V}(n)$ of finite type for all $n\ge0$. If\[
\alpha_{i}:C_{i}\to P_{\mathcal{V}}(\forgetful{C_{i}})\]
are the classifying maps with\[
\hat{\alpha}:\ilimit\forgetful{C_{i}}\to\ilimit\forgetful{P_{\mathcal{V}}(\forgetful{C_{i}})}\]
the induced map, and if \[
\mu:\forgetful{P_{\mathcal{V}}(\ilimit\forgetful{C_{i}})}\to\ilimit\forgetful{P_{\mathcal{V}}(\forgetful{C_{i}})}\]
is the inclusion defined in proposition~\ref{pro:injectioninvercofree},
then\[
\mu\left(\forgetful{\ilimit C_{i}}\right)=\mu\left(\forgetful{P_{\mathcal{V}}(\ilimit\forgetful{C_{i}})}\right)\cap\hat{\alpha}\left(\ilimit\forgetful{C_{i}}\right)\subseteq\ilimit{\forgetful{P_{\mathcal{V}}(\forgetful{C_{i}})}}\]

A corresponding results holds in the category $\icoalgcat$ after
consistently replacing the functor $P_{\mathcal{V}}(*)$ by $L_{\mathcal{V}}(*)$.\end{cor}
\begin{rem*}
The naive way to construct $\ilimit C_{i}$ is to try to equip $\ilimit\forgetful{C_{i}}$
with a coproduct --- a process that fails because we only get a map\[
\ilimit\forgetful{C_{i}}\to\prod_{n\ge0}\homzs n(\mathcal{V}(n),\ilimit(C_{i}^{\otimes n}))\ne\prod_{n\ge0}\homzs n(\mathcal{V}(n),(\ilimit C_{i})^{\otimes n})\]
which is not a true coalgebra structure.

Corollary\,\ref{cor:inverselimitintersectiontwo} implies that this
naive procedure \emph{almost} works. Its failure is precisely captured
by the degree to which\[
\forgetful{P_{\mathcal{V}}(\ilimit\forgetful{C_{i}})}\ne\ilimit{\forgetful{P_{\mathcal{V}}(\forgetful{C_{i}})}}\]
\end{rem*}
\begin{proof}
This follows immediately from the exact sequence~\ref{eq:sixterminvlimexact}.
\end{proof}
Our main result
\begin{thm}
\label{thm:aequivforgetfuldirectsummand}Let $\{f_{i}\}:\{A\}\to\{C_{i}\}$
be a morphism in $\iircoalgcat$ over a $\Sigma$-cofibrant operad
$\mathcal{V}=\{\mathcal{V}(n)\}$ with $\mathcal{V}(n)$ of finite
type for all $n\ge0$. Let$ $
\begin{enumerate}
\item $\{A\}$ be the constant object
\item the $\{f_{i}\}$ be cofibrations for all $i$
\end{enumerate}
Then $\{f_{i}\}$ induces an inclusion $f=\ilimit\{f_{i}\}:A\to\ilimit\{C_{i}\}$
and the sequence \[
0\to\forgetful A\xrightarrow{\ilimit f_{i}}\forgetful{\ilimit C_{i}}\to\forgetful{\ilimit(C_{i}/A)}\to0\]
is exact. In particular, if $ $$\forgetful{\ilimit(C_{i}/A)}$ is
contractible, then $\ilimit f_{i}$ is a weak equivalence.\end{thm}
\begin{proof}
We will consider the case of $\coalgcat$ --- the other case follows
by a similar argument. 

The inclusion\[
\forgetful{\ilimit C_{i}}\subseteq\ilimit\forgetful{C_{i}}\]
from corollary~\ref{cor:intersectionalreadycoalg}, and the left-exactness
of filtered limits in $\chaincat$ implies the left-exactness of the
filtered limits in $\iircoalgcat$, and that the inclusion\[
\forgetful A\hookrightarrow\forgetful{\ilimit C_{i}}\]
is a cofibration in $\chaincat$.

The fact that\[
\forgetful{\ilimit C_{i}}=\forgetful{\bigcap_{i=1}^{\infty}q_{i}^{-1}(\alpha_{i}(C_{i}))}\subseteq\forgetful{L_{\mathcal{V}}(\ilimit\forgetful{C_{i}})}\]
from the same corollary and the diagram \[
\xymatrix{{q_{j}^{-1}(\alpha_{j}(C_{j}))}\ar@{^{(}->}[d]\ar[r]^{h\quad} & {u_{j}^{-1}(\alpha'_{j}(C_{j}/A))}\ar@{^{(}->}[d]\\
{L_{\mathcal{V}}(\ilimit\forgetful{C_{i}})}\ar@{->>}[r]\ar[d]_{q_{j}} & {L_{\mathcal{V}}(\ilimit\forgetful{C_{i}/A})}\ar[d]^{u_{j}}\\
{L_{\mathcal{V}}(\forgetful{C_{j}})}\ar@{->>}[r] & {L_{\mathcal{V}}(\forgetful{C_{j}/A})}\\
{C_{j}}\ar[u]^{\alpha_{j}}\ar@{->>}[r] & {C_{j}/A}\ar[u]_{\alpha'_{j}}}
\]
shows that the map $h$ is surjective. The conclusion follows. 
\end{proof}
\bibliographystyle{amsplain}

\providecommand{\bysame}{\leavevmode\hbox to3em{\hrulefill}\thinspace}
\providecommand{\MR}{\relax\ifhmode\unskip\space\fi MR }
\providecommand{\MRhref}[2]{%
  \href{http://www.ams.org/mathscinet-getitem?mr=#1}{#2}
}
\providecommand{\href}[2]{#2}

\end{document}